\documentclass[11pt]{article}

\usepackage{amssymb,amsbsy,amsthm,amsmath,graphicx,mathrsfs,epsfig,graphics,times}

\input xy

\xyoption{all}

\begin{document}

\title{\textbf{Pleasant extensions retaining algebraic
structure, I}}
\author{Tim Austin}
\date{}

\maketitle


\newenvironment{nmath}{\begin{center}\begin{math}}{\end{math}\end{center}}

\newtheorem{thm}{Theorem}[section]
\newtheorem*{thm*}{Theorem}
\newtheorem{lem}[thm]{Lemma}
\newtheorem{prop}[thm]{Proposition}
\newtheorem{cor}[thm]{Corollary}
\newtheorem{conj}[thm]{Conjecture}
\newtheorem{dfn}[thm]{Definition}
\newtheorem{prob}[thm]{Problem}
\newtheorem{ques}[thm]{Question}
\theoremstyle{remark}
\newtheorem*{ex}{Example}


\newcommand{\s}{\sigma}
\renewcommand{\O}{\Omega}
\renewcommand{\S}{\Sigma}
\newcommand{\co}{\mathrm{co}}
\newcommand{\e}{\mathrm{e}}
\newcommand{\eps}{\varepsilon}
\renewcommand{\d}{\mathrm{d}}
\newcommand{\im}{\mathrm{i}}
\renewcommand{\k}{\kappa}
\renewcommand{\l}{\lambda}
\newcommand{\G}{\Gamma}
\newcommand{\g}{\gamma}
\renewcommand{\L}{\Lambda}
\renewcommand{\a}{\alpha}
\renewcommand{\b}{\beta}
\newcommand{\Sone}{\mathrm{S}^1}

\newcommand{\Aut}{\mathrm{Aut}}
\renewcommand{\Pr}{\mathrm{Pr}}
\newcommand{\Hom}{\mathrm{Hom}}
\newcommand{\id}{\mathrm{id}}
\newcommand{\BSp}{\mathsf{BSp}}
\newcommand{\WBSp}{\mathsf{WBSp}}
\newcommand{\Sys}{\underline{\mathsf{Sys}}}
\newcommand{\Rep}{\underline{\mathsf{Rep}}}
\newcommand{\Lat}{\mathrm{Lat}}
\newcommand{\FLat}{\mathrm{FLat}}
\newcommand{\SG}{\mathrm{SG}}
\newcommand{\CSG}{\mathrm{CSG}}
\newcommand{\Clos}{\mathrm{Clos}}
\newcommand{\pro}{\mathrm{pro}}
\newcommand{\nil}{\mathrm{nil}}
\newcommand{\rat}{\mathrm{rat}}
\newcommand{\Ab}{\mathrm{Ab}}
\newcommand{\CIS}{\mathrm{CIS}}

\newcommand{\bbN}{\mathbb{N}}
\newcommand{\bbR}{\mathbb{R}}
\newcommand{\bbZ}{\mathbb{Z}}
\newcommand{\bbQ}{\mathbb{Q}}
\newcommand{\bbT}{\mathbb{T}}
\newcommand{\bbD}{\mathbb{D}}
\newcommand{\bbC}{\mathbb{C}}
\newcommand{\bbP}{\mathbb{P}}

\newcommand{\A}{\mathcal{A}}
\newcommand{\B}{\mathcal{B}}
\newcommand{\C}{\mathcal{C}}
\newcommand{\E}{\mathcal{E}}
\newcommand{\F}{\mathcal{F}}
\newcommand{\I}{\mathcal{I}}
\newcommand{\calL}{\mathcal{L}}
\renewcommand{\P}{\mathcal{P}}
\newcommand{\U}{\mathcal{U}}
\newcommand{\W}{\mathcal{W}}
\newcommand{\Y}{\mathcal{Y}}
\newcommand{\Z}{\mathcal{Z}}

\newcommand{\frH}{\mathfrak{H}}
\newcommand{\frM}{\mathfrak{M}}

\newcommand{\bfX}{\mathbf{X}}
\newcommand{\bfY}{\mathbf{Y}}
\newcommand{\bfZ}{\mathbf{Z}}
\newcommand{\bfW}{\mathbf{W}}
\newcommand{\bfV}{\mathbf{V}}

\newcommand{\rmC}{\mathrm{C}}
\newcommand{\rmH}{\mathrm{H}}
\newcommand{\T}{\mathrm{T}}

\newcommand{\sfC}{\mathsf{C}}
\newcommand{\sfE}{\mathsf{E}}
\newcommand{\sfF}{\mathsf{F}}
\newcommand{\sfV}{\mathsf{V}}
\newcommand{\sfW}{\mathsf{W}}
\newcommand{\sfY}{\mathsf{Y}}
\newcommand{\sfZ}{\mathsf{Z}}

\newcommand{\uhr}{\!\!\upharpoonright}
\newcommand{\into}{\hookrightarrow}
\newcommand{\onto}{\twoheadrightarrow}

\newcommand{\bb}[1]{\mathbb{#1}}
\newcommand{\bs}[1]{\boldsymbol{#1}}
\newcommand{\fr}[1]{\mathfrak{#1}}
\renewcommand{\bf}[1]{\mathbf{#1}}
\renewcommand{\sf}[1]{\mathsf{#1}}
\renewcommand{\rm}[1]{\mathrm{#1}}
\renewcommand{\cal}[1]{\mathcal{#1}}

\renewcommand{\t}[1]{\tilde{#1}}

\newcommand{\mvee}{\hbox{$\bigvee$}}

\newcommand{\lfr}{\lfloor}
\newcommand{\rfr}{\rfloor}

\newcommand{\fin}{\nolinebreak\hspace{\stretch{1}}$\lhd$}
\newcommand{\tick}[1]{\nolinebreak\hspace{\stretch{1}}$\surd_{\mathrm{#1}}$}

\begin{abstract}
In the recent papers~\cite{Aus--nonconv}
and~\cite{Aus--newmultiSzem} we introduced some new techniques for
constructing an extension of a probability-preserving system $T:\bbZ^d\curvearrowright
(X,\mu)$ that enjoys certain desirable properties
in connexion with the asymptotic behaviour of some related nonconventional ergodic
averages.

The present paper is the first of two that will explore various
refinements and extensions of these ideas.  This first part is
dedicated to some much more general machinery for the construction
of extensions that can be used to recover the results
of~\cite{Aus--nonconv,Aus--newmultiSzem}.  It also contains two
relatively simple new applications of this machinery to the study of
certain families of nonconventional averages, one in discrete and
one in continuous time (convergence being a new result for the
latter).

In the forthcoming second part~\cite{Aus--lindeppleasant2} we will
introduce the problem of describing the characteristic factors and
the limit of the linear nonconventional averages
\[\frac{1}{N}\sum_{n=1}^N \prod_{i=1}^kf_i\circ T^{n\bf{p}_i}\]
when the directions $\bf{p}_1$, $\bf{p}_2$, \ldots, $\bf{p}_k \in
\bbZ^d$ are not assumed to be linearly independent, and provide a
fairly detailed solution in the case when $k = 3$, $d= 2$ and any
pair of directions is linearly independent.  This will then be used to prove the convergence in $L^2(\mu)$ of the quadratic nonconventional averages
\[\frac{1}{N}\sum_{n=1}^N (f_1\circ T_1^{n^2})(f_2\circ T_1^{n^2}T_2^n).\]
\end{abstract}

\parskip 0pt

\tableofcontents

\parskip 7pt

\section{Introduction}

Suppose that $T_1$, $T_2$, \ldots, $T_k \curvearrowright (X,\mu)$ is
a system of commuting invertible transformations on a standard Borel
probability space.  To such a system we can associate various
`nonconventional' ergodic averages, such as the averages
\[\frac{1}{N}\sum_{n=1}^N\prod_{i=1}^k f_i\circ T_i^n\]
or their more complicated relatives of the form
\[\frac{1}{N}\sum_{n=1}^N\prod_{i=1}^k f_i\circ T^{p_i(n)}\]
for an action $T:\bbZ^d\curvearrowright (X,\mu)$ and polynomial
mappings $p_i:\bbZ\longrightarrow \bbZ^d$ for $i=1,2,\ldots,k$
(we sometimes refer to these as `linear' and `polynomial' averages respectively).

The linear averages have been the subject of considerable study
since they first emerged from Furstenberg's ergodic theoretic
approach to Szemer\'edi's Theorem and
generalizations~\cite{Fur77,FurKat78}, and more recently their
polynomial relatives have received similar attention since the
extension of Furstenberg's work to an appropriate polynomial
setting~\cite{Ber87,BerLei96}.  These applications to Arithmetic
Ramsey Theory typically require that certain related scalar averages
`stay large' as $N\to \infty$, but it quickly became clear that the
more fundamental question of their norm convergence in $L^2(\mu)$
posed an interesting challenge in its own right. After several
important partial
results~\cite{ConLes84,ConLes88.1,ConLes88.2,Mei90,Zha96,FurWei96,HosKra01,HosKra05,Zie05},
the convergence of the linear averages in general was settled by Tao
in~\cite{Tao08(nonconv)}.  By contrast, our understanding of the
polynomial case remains poor.

In~\cite{Aus--nonconv} we gave a new proof of convergence in the
linear case, more classically ergodic theoretic than Tao's (which
relies on a conversion of the problem into an equivalent
quantitative assertion about the shift transformations on
$(\bbZ/N\bbZ)^d$).  In this paper and its 
sequel~\cite{Aus--lindeppleasant2} we shall
further develop the methods of~\cite{Aus--nonconv} to provide some
more versatile machinery, and illustrate its use with a proof of
norm convergence for the new polynomial instance
\[\frac{1}{N}\sum_{n=1}^N(f_1\circ T_1^{n^2})(f_2\circ T_1^{n^2}T_2^n).\]

Most analyses of such questions rest on the fundamental notion of
`characteristic factors' for a system of averages, first made
explicit by Furstenberg and Weiss in their work~\cite{FurWei96} on a
case of polynomial averages involving only a single underlying
transformation.  Here a tuple of factors $\xi_i:(X,\mu,T)\to
(Y_i,\nu_i,S_i)$, $i=1,2,\ldots,k$, of a $\bbZ^d$-system will be
termed \textbf{characteristic} for some tuple of polynomial mappings
$p_i:\bbZ\longrightarrow \bbZ^d$ if
\[\frac{1}{N}\sum_{n=1}^N\prod_{i=1}^df_i \circ T^{p_i(n)}\sim \frac{1}{N}\sum_{n=1}^N\prod_{i=1}^d\sfE_\mu(f_i\,|\,\xi_i) \circ T^{p_i(n)},\]
for any $f_1$, $f_2$, \ldots, $f_d \in L^\infty(\mu)$, where we
write $f_N \sim g_N$ to denote that $\|f_N - g_N\|_2 \to 0$ as
$N\to\infty$.

In the case $d = 1$ (so we deal only with powers of a single
transformation) the possible structures of such characteristic
factors have been completely understood as those of pro-nilsystems,
first in the case of linear averages~\cite{HosKra05,Zie07} and then
also for their higher-degree polynomial
relatives~\cite{HosKra05poly,Lei05}.

However, for general $d$, even
when each $p_i$ is a linear mapping the possible characteristic
factors seem much more complicated; although it can be shown
abstractly that minimal characteristic factors exist, they have so
far largely resisted useful description.  The chief innovations of the
papers~\cite{Aus--nonconv,Aus--newmultiSzem} were methods for
constructing an extension a system $T:\bbZ^d\curvearrowright
(X,\mu)$ to a larger system $\t{T}:\bbZ^d\curvearrowright
(\t{X},\t{\mu})$ (which depends on the linear mappings $p_i$) in
which a characteristic tuple of factors for the convergence of these
linear averages could be found with an especially simple structure.
Knowing this structure then enabled (together with an appeal to
various existing machinery) new ergodic-theoretic proofs of first
the convergence of such averages in $L^2(\mu)$ and then the
Multidimensional Multiple Recurrence Theorem.

Effectively, this approach shifts our viewpoint from the setting of
an individual system of commuting transformations (and factors
thereof) to the category of all such systems. While a given system
can fail to exhibit among its own factors all of the structural
features that can be exploited to examine nonconventional averages,
these emerge upon passing to a sufficiently enriched extension of
that system (or, equivalently, upon considering more general
joinings of it to other systems of various special kinds). These extensions with improved
behaviour do not seem to be canonical: the constructions below proceed in several steps, involving both some arbitrary choices and a somewhat arbitrary order.  In this
sequence of papers we will generally refer to extensions that admit
some simple characteristic factors for a given system of averages as
\textbf{pleasant} extensions of the original system, although the
particular class of `simple' characteristic factors that we use will
vary from one context to the next.

In fact, the strategy of passing to an extension where the behaviour
of nonconventional averages is more easily described also has a
precedent in various earlier papers, notably the
work~\cite{FurWei96} of Furstenberg and Weiss (we will return to the
relationship between their work and ours several times later).
However, while in their work it can later be proved that the
characteristic factors of the original system must have taken the
same form as those obtained in the extension, in the multidimensional
setting the passage to an extension leads to a genuine reduction in
complexity of description of the characteristic factors, and recent
works in this area have begun to exploit this idea much more
extensively. (For example, Bernard Host~\cite{Hos09} has given a new
and rather more efficient construction of an extension with much the
same desirable properties as those constructed
in~\cite{Aus--nonconv}.)

Here we develop further this approach to nonconventional averages.
An undesirable feature of the constructions of
extensions in~\cite{Aus--nonconv,Aus--newmultiSzem,Hos09} is that
the lifted commuting transformations $\t{T}_i$ generally lose any
algebraic structure that might have been known to hold a priori
among the $T_i$.  Two particular kinds of structure that can be of
interest for applications are:
\begin{itemize}
\item the existence of roots, such as some $S_i$ such that $T_i =
S_i^2$;
\item linear relations that may hold among the original
transformations $T_i$, such as $T_3 = T_1T_2$.
\end{itemize}

In this paper and its sequels we will begin to see what pleasant
extensions can be found that retain such additional features. Concerning the existence of roots, in Section~\ref{subs:roots} we will find that one can recover essentially the same result as in~\cite{Aus--nonconv}:

\begin{thm}[Pleasant extensions of linearly independent linear averages]\label{thm:char-retaining-roots}
Any $\bbZ^d$-system $(X,\mu,T)$ has an extension
$\pi:(\t{X},\t{\mu},\t{T})\to (X,\mu,T)$ such that for any linearly
independent $\bf{p}_1$, $\bf{p}_2$, \ldots, $\bf{p}_k\in \bbZ^d$ the
averages
\[\frac{1}{N}\sum_{n=1}^N\prod_{i=1}^d f_i\circ \t{T}^{n\bf{p}_i},\quad\quad f_1,f_2,\ldots,f_d \in L^\infty(\t{\mu}),\]
admit a characteristic tuple of factors of the form
\[\xi_i :=\zeta_0^{\t{T}^{\bf{p}_i}}\vee\bigvee_{j\neq
i}\zeta_0^{\t{T}^{\bf{p}_i - \bf{p}_j}},\quad\quad i=1,2,\ldots,d,\] where for a
transformation $S$ we write $\zeta_0^S$ for some factor map
generating, up to negligible sets, the $\s$-algebra of sets left
invariant by $S$.
\end{thm}

(The concepts and notation appearing in this statement will be introduced more carefully in Section~\ref{sec:background}.)

The point here is that, setting $\L:= \bbZ\bf{p}_1 + \cdots +
\bbZ\bf{p}_d \leq \bbZ^d$, the older constructions give only an
extension of the subaction $T^{\ \uhr\L}$.  Notice that we obtain a
single extended system $(\t{X},\t{\mu},\t{T})$ that enjoys the above
simplified characteristic factors for every tuple of linearly
independent directions.  We will refer to such an extension as a
\textbf{pleasant extension for linearly independent linear
averages}.

More importantly than this result, this first paper will introduce
various general ideas needed in preparation for the more
sophisticated results of the sequel~\cite{Aus--lindeppleasant2}, especially the notion of
`satedness' for probability-preserving systems (Subsection 3.1
below) which will be relied on repeatedly in proving all the main
results of that paper.

In~\cite{Aus--lindeppleasant2} we shall begin to address the second
kind of algebraic structure listed above. We will examine one
simple case in detail, but conjecture that our methods could
eventually be extended to a much more general result. We will
consider the case of three directions $T^{\bf{p}_1}$,
$T^{\bf{p}_2}$, $T^{\bf{p}_3}$ in a $\bbZ^2$-system that are in
general position with the origin $\bs{0}$: that is, such that no
three of the points $\bf{p}_1$, $\bf{p}_2$, $\bf{p}_3$ and $\bs{0}$
lie on a line. For the associated linear averages we will show how
to construct an extended $\bbZ^2$-system in which the characteristic
factors take a special form, which will give us a notion of pleasant
extensions for triple linear averages subject to this kind of linear
dependence.

Of course, we do pay a price for insisting that the
$\bbZ^2$-structure of the action be preserved, in that the
characteristic factors we eventually obtain are not as simple as the
pure joins of isotropy factors that emerge in the linearly
independent case.  The additional ingredients we need are $\bbZ^2$-actions given by pairs of commuting rotations on two-step pro-nilmanifolds (or, to be precise, direct integrals of such actions).

\begin{thm}[Pleasant extensions for linearly dependant triple linear averages]\label{thm:char-three-lines-in-2D}
Any $\bbZ^2$-system $(X,\mu,T)$ has an extension
$\pi:(\t{X},\t{\mu},\t{T}) \to (X,\mu,T)$ such that for any
$\bf{p}_1$, $\bf{p}_2$, $\bf{p}_3 \in \bbZ^2$ that are in general
position with the origin the averages
\[\frac{1}{N}\sum_{n=1}^N(f_1\circ \t{T}^{n\bf{p}_1})(f_2\circ \t{T}^{n\bf{p}_2})(f_3\circ \t{T}^{n\bf{p}_3}),\quad\quad f_1,f_2,f_3 \in L^\infty(\t{\mu}),\]
admit a characteristic triple of factors of the form
\[\xi_i = \zeta_0^{\t{T}^{\bf{p}_i}}\vee \zeta_0^{\t{T}^{\bf{p}_i - \bf{p}_j}}\vee\zeta_0^{\t{T}^{\bf{p}_i - \bf{p}_k}}\vee \zeta_{\nil,2}^{\t{T}},\quad\quad i=1,2,3,\]
where $\{i,j,k\} = \{1,2,3\}$ and $\zeta_{\nil,2}^{\t{T}}$ is the maximal factor of $(X,\mu,T)$ generated by direct integrals of two-step nilsystems (the precise meaning of this will be elaborated in~\cite{Aus--lindeppleasant2}).
\end{thm}

From this we will also be able to deduce a pleasant-extensions result for certain double quadratic averages, following an application of the well-known van der Corput estimate.

\begin{thm}[Pleasant extensions for some double quadratic averages]\label{thm:char-poly}
Any system of two commuting transformations $T_1,T_2\curvearrowright
(X,\mu)$ has an extension $\pi:(\t{X},\t{\mu},\t{T}_1,\t{T}_2) \to
(X,\mu,T_1,T_2)$ in which the averages
\[\frac{1}{N}\sum_{n=1}^N(f_1\circ \t{T}_1^{n^2})(f_2\circ \t{T}_1^{n^2}\t{T}_2^n),\quad\quad f_1,f_2 \in L^\infty(\t{\mu}),\]
admit characteristic factors of the form
\[\xi_1 = \xi_2 := \bigvee_{m\geq 1}\zeta_0^{\t{T}^m_1}\vee \zeta_0^{\t{T}_2}\vee \zeta_{\nil,2}^{\t{T}}.\]
\end{thm}

This, in turn, gives us enough control over these quadratic averages to complete a proof of their norm convergence.

\begin{thm}\label{thm:polyconv}
If $T_1,T_2:\bbZ\curvearrowright (X,\mu)$ commute then the averages
\[\frac{1}{N}\sum_{n=1}^N (f_1\circ T_1^{n^2})(f_2\circ T_1^{n^2}T_2^n)\]
converge in $L^2(\mu)$ as $N\to\infty$ for any $f_1,f_2 \in
L^\infty(\mu)$.
\end{thm}

Although this new convergence result is modest in itself --- it is
only one special case of the much more general conjecture of norm
convergence for all multidimensional polynomial averages, which
remains out of reach for the time being --- we suspect that the
methods developed in this paper and its sequel will ultimately have more far-reaching relevance to this question, and potentially to other
questions on the structure of joinings between different classes of
system in the ergodic theory of $\bbZ^d$-actions.

\textbf{Remark}\quad A long time passed between the submission of the present paper and its acceptance.  During that interval, the author found some other applications of the notion of satedness to the study of nonconventional averages, including to some actions of non-discrete or non-Abelian groups:~\cite{Aus--normconvctstime,Aus--joiningequidistnilpotent,Aus--commute}. \fin

\textbf{Acknowledgements}\quad My thanks go to Vitaly Bergelson,
Zolt\'an Buczolich, Bernard Host, Bryna Kra, Mariusz Lema\'nczyk,
Emmanuel Lesigne, Thierry de la Rue, Terence Tao, Jean-Paul Thouvenot, Dave Witte
Morris and Tamar Ziegler for several helpful discussions, and to the
Mathematical Sciences Research Institute (Berkeley) for its
hospitality during the 2008 program on Ergodic Theory and Additive
Combinatorics.  My thanks go also to the referee, who suggested a number of improvements.

\section{Background definitions and general results}\label{sec:background}

\subsection{Measurable functions and probability kernels}

We will work exclusively in the category of standard Borel
probability spaces $(X,\S_X,\mu)$, and so will often suppress
mention of their $\s$-algebras.  Given standard Borel spaces $X$,
$Y$ and $Z$, a completed Borel probability measure $\mu$ on $X$, a
$\s$-subalgebra $\Xi\subseteq \S_X$ and measurable functions
$\phi:X\to Y$, $\psi:X\to Z$, we will write that $\psi$ is
\textbf{$\mu$-virtually $\Xi$-measurable} if there is some
$\Xi$-measurable map $\psi_1:X\to Z$ such that $\psi= \psi_1$
$\mu$-almost everywhere, or similarly that $\psi$ is
\textbf{$\mu$-virtually a function of $\phi$} if there is some
measurable function $\theta:Y \to Z$ such that $\psi =
\theta\circ\phi$ $\mu$-almost everywhere; these two definitions are related by the usual
correspondence (up to negligible sets) between $\s$-subalgebras and
maps that generate them for standard Borel spaces.

Factor maps from one probability space to another comprise the
simplest class of morphisms between such spaces, but we will
sometimes find ourselves handling also a weaker class of morphisms.
Suppose that $Y$ and $X$ are standard Borel spaces.  Then by a
\textbf{probability kernel from $Y$ to $X$} we understand a function
$P:Y\times\S_X \to [0,1]$ such that
\begin{itemize}
\item the map $y\mapsto P(y,A)$ is $\S_Y$-measurable for every $A \in
\S_X$;
\item the map $A \mapsto P(y,A)$ is a probability measure on $\S_X$
for every $y \in Y$.
\end{itemize}
The first of the above conditions is then the natural sense in which
the assignment $y \mapsto P(y,\,\cdot\,)$ of a probability
distribution is measurable in $y$; indeed, a popular alternative
definition of probability kernel is as a measurable function from
$Y$ to the set $\Pr\,X$ of Borel probability measures on $X$.  In
ergodic theory this notion is also often referred to as a
`quasifactor': see, for example, Chapter 8 of Glasner~\cite{Gla03},
where this alternative convention and notation are used. We will
write $P:Y \stackrel{\rm{p}}{\to} X$ when $P$ is a probability
kernel from $Y$ to $X$.

Given a kernel $P:Y \stackrel{\rm{p}}{\to} X$ and a probability
measure $\nu$ on $Y$, we define the \textbf{pushforward} measure
$P_\# \nu$ on $X$ by
\[P_\# \nu(A) := \int_YP(y,A)\,\nu(\d y);\]
this measure on $X$ can be interpreted as the law of a member of $X$
selected randomly by first selecting a member of $Y$ with law $\nu$
and then selecting a member of $X$ with law $P(y,\,\cdot\,)$. This
extends standard deterministic notation: given a measurable function
$\phi:Y \to X$, we may associate to it the deterministic probability
kernel given by $P(y,\,\cdot\,) = \delta_{\phi(y)}$ (the point mass
at the image of $y$ under $\phi$), and now $P_\#\nu$ is the usual
push-forward measure $\phi_\#\nu$.

Certain special probability kernels naturally serve as adjoints to
probability-preserving maps, in the sense of the following theorem.

\begin{thm}
Suppose that $Y$ and $X$ are standard Borel spaces, that $\mu$ is a
probability measure on $X$ and that $\phi:X \to Y$ is a measurable
map. Then, denoting the pushforward $\phi_\# \mu$ by $\nu$, there is
a $\nu$-almost surely unique probability kernel $P:Y
\stackrel{\rm{p}}{\to} X$ such that $\mu = P_\#\nu$ and which
\textbf{represents the conditional expectation with respect to
$\phi$}: for any $f \in L^1(\mu)$, the function
\[x_1 \mapsto \int_Xf(x)\,P(\phi(x_1),\d x)\]
is a version of the $\mu$-conditional expectation of $f$ with
respect to $\phi^{-1}(\S_Y)$.

We also write that this $P$ \textbf{represents the disintegration of
$\mu$ over $\phi$}.  A general probability kernel
$P:Y\stackrel{\rm{p}}{\to} X$ represents the disintegration over
$\phi$ of some measure that pushes forward onto $\nu$ if and only if
$\int_AP(x,\,\cdot\,)\,\nu(\d y)$ and $\int_BP(y,\,\cdot\,)\,\nu(\d
y)$ are mutually singular whenever $A\cap B = \emptyset$.
\end{thm}

\textbf{Proof}\hspace{5pt} See Theorem 6.3 in
Kallenberg~\cite{Kal02}. \qed

\subsection{Systems, subactions and factors}

In this paper we shall spend a great deal of time passing up and
down from systems to extensions or factors.  Moreover, sometimes one
system will appear as a factor of a `larger' system in \emph{several
different ways} (most obviously, when we work with a system
recovered under the different coordinate projections from some
self-joining).  For this reason the notational abuse of referring to
one system as a factor of another but leaving the relevant factor
map to the understanding of the reader, although popular and useful
in modern ergodic theory, seems dangerous here, and we shall
carefully avoid it.  In its place we substitute the alternative
abuse, slightly safer in our circumstances, of often referring only
to the factor maps we use, and leaving either their domain or target
systems to the reader's understanding.  Let us first set up some
notation to support this practice.

If $\G$ is a l.c.s.c. group, then by a \textbf{$\G$-system} (or, if $\G$
is clear, just a \textbf{system}) we understand a
probability-preserving action
$T:\G\curvearrowright (X,\mu)$ on a standard Borel probability
space that is jointly measurable (that is, it is measurable as a map $\G\times X\to X$).  We will often alternatively denote this space and action by
$(X,\mu,T)$, or by a single boldface letter such as $\bfX$. If
$\L\leq \G$ is a closed subgroup we denote by $T^{\
\uhr\L}:\L\curvearrowright (X,\mu)$ the action defined by $(T^{\
\uhr\L})^\g := T^\g$ for $\g \in \L$, and refer to this as a
\textbf{subaction}, and if $\bfX = (X,\mu,T)$ is a $\G$-system
we write similarly $\bfX^{\ \uhr\L}$ for the system $(X,\mu,T^{\
\uhr\L})$ and refer to it as a \textbf{subaction system}.

A \textbf{factor} from one system $(X,\mu,T)$ to another $(Y,\nu,S)$
is a Borel map $\pi:X\to Y$ such that $\pi_\#\mu = \nu$ and
$\pi\circ T^\g = S^\g\circ \pi$ for all $\g \in \G$. Given such a
factor, we sometimes write $T|_{\pi}$ to denote the action $S$ with
which $T$ is intertwined by $\pi$.

Any factor $\pi:\bfX\to \bfY$ specifies a globally $T$-invariant
$\s$-subalgebra of $\S_X$ in the form of $\pi^{-1}(\S_Y)$. Two
factors $\pi$ and $\psi$ are \textbf{equivalent} if these
$\s$-subalgebras of $\S_X$ that they generate are equal up to
$\mu$-negligible sets, in which case we shall write $\pi \simeq
\psi$; this clearly defines an equivalence relation among factors.

It is a standard fact that in the category of standard Borel spaces
equivalence classes of factors are in bijective correspondence with
equivalence classes of globally invariant $\s$-subalgebras under the
relation of equality modulo negligible sets. A treatment of these
classical issues may be found, for example, in Chapter 2 of
Glasner~\cite{Gla03}.  Given a globally invariant $\s$-subaglebra in
$\bfX$, a choice of factor $\pi:\bfX\to \bfY$ generating that
$\s$-subalgebra will sometimes be referred to as a
\textbf{coordinatization} of the $\s$-subalgebra.  Importantly for
us, some choices of coordinatizing factor $\pi$ may reveal some
additional structure of the factors more clearly than others. For
this reason, given one coordinatization $\pi:\bfX \to \bfY$ and an
isomorphism $\psi: \bfY \to \bfY'$, we shall sometimes refer to the
composition $\psi\circ\pi:\bfX\to \bfY'$ as a
\textbf{recoordinatization} of $\pi$.

More generally, the factor $\pi:(X,\mu,T) \to (Y,\nu,S)$
\textbf{contains} $\psi:(X,\mu,T)\to (Z,\theta,R)$ if
$\pi^{-1}(\S_Y) \supseteq \psi^{-1}(\S_Z)$ up to $\mu$-negligible
sets.  It is a classical fact that in the category of standard Borel
spaces this inclusion is equivalent to the existence of a
\textbf{factorizing} factor map $\phi:(Y,\nu,S) \to (Z,\theta,R)$
with $\psi = \phi\circ\pi$ $\mu$-a.s., and that a measurable analog
of the Schroeder-Bernstein Theorem holds: $\pi \simeq \psi$ if and
only if a single such $\phi$ may be chosen that is invertible away
from some negligible subsets of the domain and target. If $\pi$
contains $\psi$ we shall write $\pi \succsim \psi$ or $\psi \precsim
\pi$, and shall write $\psi|_\pi$ for a choice of the factorizing map $\phi$. It is clear that (up to set-theoretic niceties) this defines a
partial order on the class of $\simeq$-equivalence classes of
factors of a given system. We will extend the above terminology to
that of \textbf{coordinatizations} and \textbf{recoordinatizations}
of families of factors of a system in terms of the appropriate
commutative diagram of isomorphisms.

If $\pi:(X,\mu,T)\to (Y,\nu,S)$ and $\psi:(X,\mu,T)\to (Z,\theta,R)$ are any two factors, then the corresponding $\s$-subalgebras $\pi^{-1}(\S_Y)$ and $\psi^{-1}(\S_Z)$ generate another globally $T$-invariant $\s$-subalgebra in $\bfX$.  It therefore corresponds to another factor of $\bfX$, which is the minimal factor of $\bfX$ subject to containing both $\pi$ and $\psi$.  We write $\pi\vee\psi$ for a choice of coordinatization for this factor: one such choice is always offered by the map $X\to Y\times Z: x\mapsto (\pi(x),\psi(x))$.  On the other hand, the intersection $\pi^{-1}(\S_Y)\cap \psi^{-1}(\S_Z)$ is also a globally $T$-invariant $\s$-subalgebra in $\bfX$: it corresponds to the maximal factor that is contained in both $\pi$ and $\psi$, and we let $\pi\wedge \psi$ denote some choice of coordinatization for it (although in this case there is no canonical such choice).

Given a $\G$-system $\bfX = (X,\mu,T)$, the $\s$-algebra $\S_X^T$ of
sets $A\in\S_X$ for which $\mu(A\triangle T^\g(A))=0$ for all $\g
\in \G$ is $T$-invariant, so defines a factor of $\bfX$. More
generally, if $\G$ is Abelian and $\L \leq\G$ is closed then we can
consider the $\s$-algebra $\S_X^{T^{\ \uhr\L}}$ generated by all
$T^{\ \uhr\L}$-invariant sets: we refer to this as the
\textbf{$\L$-isotropy factor}.  We write $\bfZ_0^{T^{\ \uhr\L}}$ for
some new system that we adopt as the target for a factor map
$\zeta_0^{T^{\ \uhr\L}}$ that generates $\S_X^{T^{\ \uhr\L}}$, and
$Z_0^{T^{\ \uhr\L}}$ for the standard Borel space underlying
$\bfZ_0^{T^{\ \uhr\L}}$.  Note that in this notation, the original system from which this factor is obtained is marked by the explicit mention of the action $T$; this should be clearer in cases in which more than one action on the same space is being considered.  Also, the condition that $\G$ be Abelian (or, more generally, that $\L\unlhd \G$) is
needed for this to be a globally $T$-invariant factor. If $T_1$ and
$T_2$ are two commuting actions of the same group $\G$ on $(X,\mu)$
then we can define a third action $T_1T_2^{-1}$ by setting
$(T_1T_2^{-1})^\g := T_1^\g T_2^{\g^{-1}}$, and in this case we
sometimes write $\zeta_0^{T_1 = T_2}:\bfX\to\bfZ_0^{T_1 = T_2}$ in
place of $\zeta_0^{T_1T_2^{-1}}:\bfX\to \bfZ_0^{T_1T_2^{-1}}$.
Similarly, if $S\subseteq \G$ and $\L$ is the group generated by
$S$, we will sometimes write $\bfZ_0^{T^{\ \uhr S}}$ for
$\bfZ_0^{T^{\ \uhr \L}}$, and similarly.

An important construction of new systems from old is that of
\textbf{relatively independent products}.  If $\bfY = (Y,\nu,S)$ is
some fixed system and $\pi_i:\bfX_i = (X_i,\mu_i,T_i)\to \bfY$ is an
extension of it for $i=1,2,\ldots,k$ then we define the relatively
independent product of the systems $\bfX_i$ over their factor maps
$\pi_i$ to be the system
\[\prod_{\{\pi_1 = \pi_2 = \ldots = \pi_k\}}\bfX_i = \Big(\prod_{\{\pi_1 = \pi_2 = \ldots = \pi_k\}}X_i,\bigotimes_{\{\pi_1 = \pi_2 = \ldots = \pi_k\}}\mu_i,T_1\times T_2\times\cdots\times T_k\Big)\]
where
\begin{multline*}
\prod_{\{\pi_1 = \pi_2 = \ldots = \pi_k\}}X_i :=
\{(x_1,x_2,\ldots,x_k)\in X_1\times X_2\times\cdots\times X_k:\\
\pi_1(x_1) = \pi_2(x_2) = \ldots = \pi_k(x_k)\},
\end{multline*}
\[\bigotimes_{\{\pi_1 = \pi_2 = \ldots = \pi_k\}}\mu_i = \int_Y \bigotimes_{i=1}^kP_i(y,\,\cdot\,)\,\nu(\d y)\]
and $P_i:Y\stackrel{\rm{p}}{\to} X_i$ is a probability kernel
representing the disintegration of $\mu_i$ over $\pi_i$. In case
$k=2$ we will write this instead as $\bfX_1\times_{\{\pi_1=
\pi_2\}}\bfX_2$, and in addition if $\bfX_1 = \bfX_2 = \bfX$ and
$\pi_1 = \pi_2 = \pi$ then we will abbreviate this further to
$\bfX\times_\pi\bfX$, and similarly for the individual spaces and
measures.

\section{Idempotent classes of systems}

\subsection{Basic properties of idempotent classes}

Central to this paper will be a systematic exploitation of a
property of certain systems according to which they can be joined to
certain other classes of system only in simple ways. This key
definition, although very abstract and very simple, will repeatedly
prove surprisingly powerful. We will introduce it after some other
preliminaries about classes of $\G$-systems.

\begin{lem}
Suppose that $\sfC$ is a class of $\G$-systems (formally, $\sfC$ is
a subcategory of $\G$-$\Sys$) that contains the trivial system and
is closed under isomorphisms, arbitrary finite joinings and inverse limits.
Then any $\G$-system $\bfX$ has an essentially unique maximal factor
in the class $\sfC$.
\end{lem}

\textbf{Proof}\quad It is clear that under the above assumption the
class
\[\{\Xi \subseteq \S_X:\ \Xi\ \hbox{is a $T$-invariant $\s$-subalgebra such that the associated factor is in $\sfC$}\}\]
is nonempty (it contains $\{\emptyset,X\}$), upwards directed
(because $\sfC$ is closed under joinings) and closed under taking
$\s$-algebra completions of increasing unions (because $\sfC$ is
closed under inverse limits).  There is therefore a maximal
$\s$-subalgebra in this set. \qed

\begin{dfn}[Idempotence]
A class of systems $\sfC$ is \textbf{idempotent} if it contains the
trivial system and is closed under isomorphisms, finite joinings and
inverse limits. In this case, we will write that $\bfX$ is a
\textbf{$\sfC$-system} if $\bfX$ is a system in the class $\sfC$,
and for arbitrary $\bfX$ we write $\zeta_\sfC^\bfX:\bfX\to\sfC\bfX$
for an arbitrarily-chosen coordinatization of its \textbf{maximal
$\sfC$-factor} given by the above lemma.

It is clear that if $\pi:\bfX\to\bfY$ then $\zeta_\sfC^\bfX \succsim
\zeta_\sfC^\bfY\circ\pi$, and so there is an essentially unique
factorizing map, which we denote by $\sfC\pi$, that makes the
following diagram commute:
\begin{center}
$\phantom{i}$\xymatrix{
&\bfX\ar[dl]_{\zeta^\bfX_\sfC}\ar[dr]^{\pi}\\
\sfC\bfX\ar[dr]_{\sfC\pi} & & \bfY\ar[dl]^{\zeta^\bfY_\sfC}\\
&\sfC\bfY.}
\end{center}

In addition, we shall abbreviate $\bfX\times_{\zeta_\sfC^\bfX} \bfX$
to $\bfX\times_\sfC \bfX$, and similarly for the individual spaces
and measures defining these systems.
\end{dfn}

The reason for this terminology lies in the observation that the
assignment $\bfX\mapsto \sfC\bfX$ defines an autofunctor of the
category $\G$-$\Sys$, and the assignment
$\bfX\mapsto\zeta_\sfC^\bfX$ defines a natural transformation from
the identity functor to this autofunctor. We can work with such
functors quite generally, and a simple definition-chase shows that a
functor $\sfF:\G\hbox{-$\Sys$}\to\G\hbox{-$\Sys$}$ together with a
natural transformation $\id_{\G\hbox{-$\Sys$}}\to \sfF$ correspond
to a class of systems $\sfC$ as above if and only if (i) $\sfF$ is
idempotent (that is, $\sfF(\sfF\bfX) = \sfF\bfX$ for every $\bfX$), (ii) the collection of systems of the form $\sfF\bfX$ is isomorphism-closed in $\G\hbox{-$\Sys$}$, and (iii) the natural transformation is the identity on the subcategory of those systems.  Indeed, it would be possible to
develop the theory of the coming sections by working entirely with
autofunctors rather than classes, but I do not know of any examples
of autofunctors that are useful as such but do not arise from
idempotent classes as above.

The name we give for our next definition is also motivated by this
relationship with functors.

\begin{dfn}[Order continuity]
A class of $\G$-systems $\sfC$ is \textbf{order continuous} if
whenever $(\bfX_{(m)})_{m\geq 0}$, $(\psi^{(m)}_{(k)})_{m\geq k\geq
0}$ is an inverse sequence of $\G$-systems with inverse limit
$\bfX$, $(\psi_{(m)})_{m\geq 0}$ we have
\[\zeta_\sfC^\bfX = \bigvee_{m\geq 0}\zeta_\sfC^{\bfX_{(m)}}\circ \psi_{(m)}:\]
that is, the maximal $\sfC$-factor of the inverse limit is simply
given by the (increasing) join of the maximal $\sfC$-factors of the
contributing systems.
\end{dfn}

\textbf{Examples}\quad The following idempotent classes will be of
particular importance (in all cases idempotence is routine to
check):
\begin{enumerate}
\item Given a fixed normal subgroup $\L\unlhd \G$, $\sfZ_0^\L$ denotes the
class of systems for which the $\L$-subaction is trivial.
\item More generally, for $\L$ as above and any $n\in \bbN$
we let $\sfZ_n^\L$ denote the class of systems on which the
$\L$-subaction is a distal tower of height at most $n$, in the sense
of direct integrals of compact homogeneous space data introduced
in~\cite{Aus--ergdirint} to allow for the case of non-ergodic
systems.
\item We can modify the previous example by placing some additional
restrictions on the permissible distal towers: for example,
$\sfZ_{\rm{Ab},n}^\L$ comprises the systems with $\L$-subaction a
distal tower of height at most $n$ and in which each isometric
extension is Abelian. We will meet other, even more restricted
idempotent classes contained within this one later.
\end{enumerate}
(We note in passing that an example of a natural class that is not
closed under joinings when $\G = \bbZ$ is the class $\sf{WM}^\perp$
of systems that are system-disjoint from all weakly-mixing systems,
as has been proved by Lema\'nczyk and Parreau in~\cite{LemPar03}.)
\fin

\textbf{Example}\quad Although all the idempotent classes that will
matter to us in this paper can be shown to be order continuous, it
may be instructive to exhibit one that is not.

Let us say that a
system $\bfX$ has a \textbf{finite-dimensional Kronecker factor} if
its Kronecker factor $\zeta_1^\bfX:X\to Z_1^\bfX$ can be
coordinatized as a direct integral (see Section 3
of~\cite{Aus--ergdirint}) of rotations on some measurably-varying
compact Abelian groups all of which can be isomorphically embedded
into a fibre repository $\bbT^D$ for some fixed $D \in \bbN$ (this
includes the possibility that the Kronecker factor is finite or
trivial).  Now let $\sfC$ be the class of $\bbZ$-systems
comprising all those that are either themselves finite-dimensional Kronecker systems, or have a Kronecker factor that is \emph{not} finite-dimensional (so we exclude just those systems that have a finite-dimensional Kronecker factor but properly contain it).  This class is idempotent.  Its closure under isomorphism and finite joinings is clear. To see that it is closed under inverse limits, suppose that $(\bfX_{(m)})_{m\geq 0}$, $(\psi^{(m)}_{(k)})_{m\geq k\geq 0}$ is an inverse sequence in $\sfC$ with inverse limit $\bfX_{(\infty)}$.  If every $\bfX_{(m)}$ is itself a Kronecker system (that is, all its ergodic components are isomorphic to rotations on compact Abelian groups), then $\bfX_{(\infty)}$ is a Kronecker system and hence lies in $\sfC$.  On the other hand, if some $\bfX_{(m)}$ is a proper extension of its Kronecker factor, then that Kronecker factor must be non-finite-dimensional since $\bfX_{(m)} \in \sfC$, and this means that $\bfX_{(\infty)}$ also has a non-finite-dimensional Kronecker factor, so also lies in $\sfC$.

However, this class $\sfC$ is not order-continuous. To see this, let $(\bfX_{(m)})_{m\geq 0}$, $(\psi^{(m)}_{(k)})_{m\geq k\geq 0}$ be an inverse sequence of ergodic rotations on finite-dimensional tori whose inverse limit $\bfX_{(\infty)}$ is an infinite-dimensional ergodic rotation, and let $\bfY$ be a nontrivial weakly mixing system.  Then $\bfX_{(\infty)} \times \bfY$ is an inverse limit of the systems $\bfX_{(m)}\times \bfY$, but this inverse limit lies in $\sfC$, whereas the maximal $\sfC$-factor of each system $\bfX_{(m)}\times \bfY$ is just $\bfX_{(m)}$, so applying $\sfC$ to the inverse sequence loses the copy of $\bfY$. \fin

\begin{dfn}[Hereditariness]
An idempotent class $\sfC$ is \textbf{hereditary} if it is also closed
under taking factors.
\end{dfn}

\textbf{Example}\quad Examples 1, 2 and 3 in the first list above
are hereditary. In the case of $\sfZ_0^\L$ this is obvious; for the
higher distal classes $\sfZ_n^\L$ or their Abelian subclasses it is
an easy consequence of the Relative Factor Structure Theorem 6.4
of~\cite{Aus--ergdirint} (applied to the relatively independent
self-joining of the $n$-step distal system in question over the
factor we wish to analyze). On the other hand, the separate example
above involving the dimensionality of the Kronecker factors is
clearly not hereditary. \fin

\begin{dfn}[Join]
If $\sfC_1$, $\sfC_2$ are idempotent classes, then the class
$\sfC_1\vee \sfC_2$ of all joinings of members of $\sfC_1$ and
$\sfC_2$ is clearly also idempotent. We call $\sfC_1\vee \sfC_2$ the
\textbf{join} of $\sfC_1$ and $\sfC_2$.
\end{dfn}

\begin{lem}[Join preserves order continuity]
If $\sfC_1$ and $\sfC_2$ are both order continuous then so is
$\sfC_1\vee \sfC_2$.
\end{lem}

\textbf{Proof}\quad Let $(\bfX_{(m)})_{m\geq 0}$,
$(\psi^{(m)}_{(k)})_{m\geq k \geq 0}$ be an inverse sequence with
inverse limit $\bfX$, $(\psi_{(m)})_{m\geq 0}$.  Then
$\zeta^\bfX_{\sfC_1\vee \sfC_2}$ is the maximal factor of $\bfX$
that is a joining of a $\sfC_1$-factor and a $\sfC_2$-factor (so, in
particular, it must be generated by its own $\sfC_1$- and
$\sfC_2$-factors), and hence it is equivalent to
$\zeta^\bfX_{\sfC_1} \vee \zeta^\bfX_{\sfC_2}$. Therefore any $f \in
L^\infty(\mu)$ that is $\zeta^\bfX_{\sfC_1} \vee
\zeta^\bfX_{\sfC_2}$-measurable can be approximated by some function
of the finite-sum form $\sum_p g_{p,1}\cdot g_{p,2}$ with each
$g_{p,i} \in L^\infty(\mu)$ being $\sfC_i$-measurable, and now since
each $\sfC_i$ is order continuous we may further approximate each
$g_{p,i}$ by some $h_{p,i} \circ \psi_{(m)}$ for a large integer $m$
and some $\sfC_i$-measurable $h_{p,i} \in L^\infty(\mu_{(m)})$.
Combining these approximations completes the proof. \qed

\textbf{Examples}\quad Of course, we can form the joins of any of
our earlier examples of idempotent classes: for example, given $\G =
\bbZ^2$ and $\bf{p}_1,\bf{p}_2,\bf{p}_3 \in \bbZ^2$ we can form
$\sfZ_0^{\bf{p}_1}\vee \sfZ_0^{\bf{p}_1-\bf{p}_2}\vee
\sfZ_0^{\bf{p}_1 - \bf{p}_3}$. This particular example and several
others like it will appear frequently throughout the rest of this
paper. We will remark shortly that joins of hereditary idempotent
classes need not be hereditary. \fin

The following terminology will also prove useful.

\begin{dfn}[Joining to an idempotent class; adjoining]
If $\bfX$ is a system and $\sfC$ is an idempotent class then a
\textbf{joining of $\bfX$ to $\sfC$} or a \textbf{$\sfC$-adjoining
of $\bfX$} is a joining of $\bfX$ and $\bfY$ for some $\bfY \in
\sfC$.
\end{dfn}

\begin{dfn}[Subjoining]
Given idempotent classes $\sfC_1$, $\sfC_2$, \ldots, $\sfC_k$, a
system $\bfX$ is a \textbf{subjoining} of $\sfC_1$, $\sfC_2$,
\ldots, $\sfC_k$ if it is a factor of a member of
$\sfC_1\vee\sfC_2\vee\cdots\vee\sfC_k$.
\end{dfn}

Finally we have reached the key definition that will drive much of the rest of this paper.

\begin{dfn}[Sated system]
Given an idempotent class $\sfC$, a system $\bfX$ is
\textbf{$\sfC$-sated} if whenever $\pi:\t{\bfX} =
(\t{X},\t{\mu},\t{T}) \to \bfX$ is an extension, the factor maps
$\pi$ and $\zeta^{\t{\bfX}}_\sfC$ are relatively
independent over $\zeta^\bfX_\sfC\circ\pi =
\sfC\pi\circ\zeta^{\t{\bfX}}_\sfC$ under $\t{\mu}$.  That is, in the commutative diagram
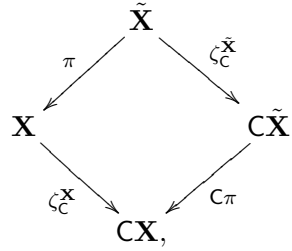
\begin{center}
$\phantom{i}$\xymatrix{
& \t{\bfX}\ar^{\zeta^{\t{\bfX}}_\sfC}[dr]\ar_\pi[dl]\\
\bfX \ar_{\zeta^\bfX_\sfC}[dr] && \sfC\t{\bfX} \ar^{\sfC\pi}[dl]\\
& \sfC\bfX, }
\end{center}
the joining of the middle two systems defined by the two factor maps from $\t{\bfX}$ must be the relatively independent product over the two maps to the common factor $\sfC\bfX$.

An inverse sequence is \textbf{$\sfC$-sated} if it has a cofinal
subsequence all of whose systems are $\sfC$-sated.
\end{dfn}

\textbf{Remark}\quad This definition has an important precedent in
Furstenberg and Weiss' notion of a `pair homomorphism' between
extensions elaborated in Section 8 of~\cite{FurWei96}. Here we shall
make much more extensive use of this basic idea. \fin

\begin{lem}
If $\sfC$ is an idempotent class, $(X,\mu,T)$ is $\sfC$-sated,
$(Y,\nu,S) \in \sfC$ and $\l$ is a $(T\times S)$-invariant
$(\mu,\nu)$-joining then $\zeta_\sfC^{(X\times Y,\l,T\times S)}
\simeq \zeta_\sfC^{(X,\mu,T)}\times \id_Y$.
\end{lem}

\textbf{Proof}\quad Let $\pi_1,\pi_2:(X\times Y,\l)\to
(X,\mu),(Y,\nu)$ be the first and second coordinate projections
respectively.  The relation $\zeta_\sfC^{(X\times Y,\l,T\times S)}
\succsim \zeta_\sfC^{(X,\mu,T)}\times \id_Y$ is clear.  On the other
hand, the $\sfC$-satedness of $(X,\mu,T)$ applied to the factor map
$\pi_1$ implies that $\pi_1$ is relatively independent from
$\zeta_\sfC^{(X\times Y,\l,T\times S)}$ over
$\zeta_\sfC^{(X,\mu,T)}\circ\pi_1$ under $\l$, and this implies the
reverse containment. This completes the proof. \qed

The crucial technical fact that turns satedness into a useful tool
is the ability to construct sated extensions of arbitrary systems.
This can be seen as a natural abstraction from Proposition 4.6
of~\cite{Aus--nonconv} and Corollary~\cite{Aus--newmultiSzem}.

\begin{thm}[Idempotent classes admit multiply sated
extensions]\label{thm:sateds-exist} If $(\sfC_i)_{i\in I}$ is a
countable family of idempotent classes then any system $\bfX_0$
admits an extension $\pi:\bfX \to \bfX_0$ such that
\begin{itemize}
\item $\bfX$ is $\sfC_i$-sated for every $i\in I$;
\item the factors $\pi$ and $\bigvee_{i\in I}\zeta^{\bfX}_{\sfC_i}$ generate the whole of
$\bfX$.
\end{itemize}
\end{thm}

We shall prove this result after a preliminary lemma.

\begin{lem}\label{lem:inverse-lim-sated}
If $\sfC$ is an idempotent class then the inverse limit of any
$\sfC$-sated inverse sequence is $\sfC$-sated.
\end{lem}

\textbf{Proof}\quad By passing to a subsequence if necessary, it
suffices to suppose that $(\bfX_{(m)})_{m\geq 0}$,
$(\psi^{(m)}_{(k)})_{m\geq k\geq 0}$ is an inverse sequence of
$\sfC$-sated systems with inverse limit $\bfX_{(\infty)}$,
$(\psi_{(m)})_{m\geq 1}$, and let $\pi:\t{\bfX} \to \bfX_{(\infty)}$
be any further extension and $f \in L^\infty(\mu_{(\infty)})$.  We
will commit the abuse of identifying such a function with its lift
to any given extension when the extension in question is obvious.
With this in mind, we need to show that
\[\sfE(f\,|\,\zeta^{\t{\bfX}}_\sfC) = \sfE(f\,|\,\zeta^{\bfX_{(\infty)}}_\sfC).\]

However, by the $\sfC$-satedness of each $\bfX_{(m)}$, we certainly
have
\[\sfE(\sfE(f\,|\,\psi_{(m)})\,|\,\zeta^{\t{\bfX}}_\sfC) = \sfE(f\,|\,\zeta^{\bfX_{(m)}}_\sfC),\]
and now as $m\to\infty$ this equation converges in $L^2(\mu)$ to
\[\sfE(f\,|\,\zeta^{\t{\bfX}}_\sfC) = \sfE\big(f\,\big|\,\lim_{m\leftarrow}\,(\zeta^{\bfX_{(m)}}_\sfC\circ\psi_{(m)})\big).\]
By monotonicity we must have
\[\zeta^{\t{\bfX}}_\sfC\succsim\zeta^{\bfX_{(\infty)}}_\sfC \succsim \lim_{m\leftarrow}\,(\zeta^{\bfX_{(m)}}_\sfC\circ \psi_{(m)}),\]
and so by sandwiching we must also have the equality of conditional
expectations desired. \qed

\textbf{Proof of Theorem~\ref{thm:sateds-exist}}\quad We first prove
this for $I$ a singleton, and then in the general case.

\quad\textbf{Step 1}\quad Suppose that $I = \{i\}$ and $\sfC_i =
\sfC$. This case will follow from a simple `energy increment'
argument.

Let $(f_r)_{r\geq 1}$ be a countable subset of the $L^\infty$-unit
ball $\{f\in L^\infty(\mu):\ \|f\|_\infty\leq 1\}$ that is dense in
this ball for the $L^2$-norm, and let $(r_i)_{i\geq 1}$ be a member
of $\bbN^\bbN$ in which every non-negative integer appears
infinitely often.

We will now construct an inverse sequence $(\bfX_{(m)})_{m\geq 0}$,
$(\psi^{(m)}_{(k)})_{m\geq k \geq 0}$ starting from $\bfX_{(0)} :=
\bfX_0$ such that each $\bfX_{(m+1)}$ is a $\sfC$-adjoining of
$\bfX_{(m)}$.  Suppose that for some $m_1 \geq 0$ we have already
obtained $(\bfX_{(m)})_{m=0}^{m_1}$, $(\psi^{(m)}_{(k)})_{m_1 \geq
m\geq k\geq 0}$ such that $\id_{X_{(m_1)}} \simeq
\zeta^{\bfX_{(m_1)}}_\sfC\vee \psi^{(m_1)}_{(0)}$.  Consider the supremum over all further extensions $\pi:\t{\bfX}\to \bfX_{(m_1)}$ of the quantity
\[\|\sfE_{\t{\mu}}(f_{r_{m_1}}\circ \psi^{(m_1)}_{(0)}\circ \pi\,|\,\zeta^{\t{\bfX}}_\sfC)\|_2^2\]
(clearly it is at most $1$), and let $\psi^{(m_1+1)}_{(m_1)}:\bfX_{(m_1 +1)}\to \bfX_{(m_1)}$ be a particular choice of extension that comes within $2^{-m_1}$ of achieving this supremum. By restricting to the possibly smaller subextension of
$\bfX_{(m_1+1)}\to\bfX_{(m_1)}$ generated by $\pi$ and
$\zeta_\sfC^{\bfX_{(m_1+1)}}$, we may assume that $\bfX_{(m_1+1)}$ is itself a
$\sfC$-adjoining of $\bfX_{(m_1)}$ and hence of $\bfX_0$. The other connecting factor maps are now determined by this $\psi^{(m_1+1)}_{(m_1)}$, so the recursion continues.

Let $\bfX_{(\infty)}$, $(\psi_{(m)})_{m \geq 0}$ be the
inverse limit of this sequence.  We have
\begin{multline*}
\id_{X_{(\infty)}} \simeq \bigvee_{m\geq 0}\psi_{(m)} \simeq
\bigvee_{m\geq 0}(\zeta_\sfC^{\bfX_{(m)}}\vee
\psi^{(m)}_{(0)})\circ\psi_{(m)}\\ \simeq \bigvee_{m\geq
0}(\zeta_\sfC^{\bfX_{(m)}}\circ\psi_{(m)})\vee\bigvee_{m\geq 0}(
\psi^{(m)}_{(0)}\circ\psi_{(m)}) \precsim
\zeta_\sfC^{\bfX_{(\infty)}}\vee \psi_{(0)},
\end{multline*}
so $\bfX_{(\infty)}$ is still a $\sfC$-adjoining of $\bfX_0$.  To
show that it is $\sfC$-sated, let $\pi:\t{\bfX}\to\bfX_{(\infty)}$
be any further extension, and suppose that $f \in
L^\infty(\mu_{(\infty)})$. We will complete the proof for Step 1 by
showing that
\[\sfE_{\t{\mu}}(f\circ\pi\,|\,\zeta_\sfC^{\t{\bfX}}) =
\sfE_{\mu_{(\infty)}}(f\,|\,\zeta_\sfC^{\bfX_{(\infty)}})\circ\pi.\]

Since $\bfX_{(\infty)}$ is a $\sfC$-adjoining of $\bfX$, this $f$
may be approximated arbitrarily well in $L^2(\mu_{(\infty)})$ by
finite sums of the form $\sum_p g_p\cdot h_p$ with $g_p$ being
bounded and $\zeta_\sfC^{\bfX_{(\infty)}}$-measurable and $h_p$
being bounded and $\psi_{(0)}$-measurable, and now by density we may
also restrict to using functions $h_p$ that are each a scalar multiple of some
$f_{r_p}\circ\psi_{(0)}$, so by continuity and multilinearity it
suffices to prove the above equality for one such product $g\cdot
(f_r\circ\psi_{(0)})$.  Since $g$ is
$\zeta_\sfC^{\bfX_{(\infty)}}$-measurable, this requirement now
reduces to
\[\sfE_{\t{\mu}}(f_r\circ\psi_{(0)}\circ\pi\,|\,\zeta_\sfC^{\t{\bfX}}) =
\sfE_{\mu_{(\infty)}}(f_r\circ\psi_{(0)}\,|\,\zeta_\sfC^{\bfX_{(\infty)}})\circ\pi.\]
Since $\zeta_\sfC^{\t{\bfX}} \succsim
\zeta_\sfC^{\bfX_{(\infty)}}\circ\pi$, this will follow if we only
show that
\[\|\sfE_{\mu_{(\infty)}}(f_r\circ\psi_{(0)}\,|\,\zeta_\sfC^{\bfX_{(\infty)}})\|_2^2 \geq \|\sfE_{\t{\mu}}(f_r\circ\psi_{(0)}\circ\pi\,|\,\zeta_\sfC^{\t{\bfX}})\|_2^2.\]

To see this, recall that $r_m = r$ for infinitely many $m$, and for all such $m$ one has
\begin{multline*}
\|\sfE_{\mu_{(m+1)}}(f_r\circ\psi^{(m)}_{(0)}\circ \psi^{(m+1)}_{(m)}\,|\,\zeta_\sfC^{\bfX_{(m+1)}})\|_2^2 \\
\geq \|\sfE_{\t{\mu}}(f_r\circ\psi^{(m)}_{(0)}\circ (\psi_{(m)}\circ\pi)\,|\,\zeta_\sfC^{\t{\bfX}})\|_2^2 - 2^{-m},
\end{multline*}
by the choice of $\psi^{(m+1)}_{(m)}:\bfX_{(m+1)}\to \bfX_{(m)}$.  Since $\|\sfE_{\mu_{(\infty)}}(f_r\circ\psi_{(0)}\,|\,\zeta_\sfC^{\bfX_{(\infty)}})\|_2^2$ is an upper bound for the left-hand value here for every $m$, we must actually have the required equality of $L^2$-norms.

\quad\textbf{Step 2}\quad The general case follows easily from Step
1 and a second inverse limit construction: choose a sequence
$(i_m)_{m\geq 1} \in I^\bbN$ in which each member of $I$ appears
infinitely often, and form an inverse sequence $(\bfX_{(m)})_{m\geq
0}$, $(\psi^{(m)}_{(k)})_{m\geq k\geq 0}$ starting from
$\bfX_{(0)}:= \bfX_0$ such that each $(\bfX_{(m)})$ is
$\sfC_{i_m}$-sated for $m\geq 1$. The inverse limit $\bfX$ is now
sated for every $\sfC_i$, by Lemma~\ref{lem:inverse-lim-sated}. \qed

\textbf{Remark}\quad Thierry de la Rue has shown me another proof of
Theorem~\ref{thm:sateds-exist} that follows very quickly from ideas
contained in his paper~\cite{LesRitdelaRue03} with Lesigne and
Rittaud, and which has now received a nice separate writeup
in~\cite{delaRue09}. The key observation is that
\begin{center}
\emph{An idempotent class $\sfC$ is hereditary if and only if
\emph{every} system is $\sfC$-sated.}
\end{center}
This in turn follows from a striking result of Lema\'nczyk, Parreau
and Thouvenot~\cite{LemParTho00} that if two systems $\bfX$ and
$\bfY$ are not disjoint then $\bfX$ shares a nontrivial factor with
the infinite Cartesian power $\bfY^{\times \infty}$.  Given now an
idempotent class $\sfC$ and a system $\bfX$, let $\sfC^\ast$ be the
hereditary idempotent class of all factors of members of $\sfC$, and
let $\bfY$ be any $\sfC$-system admitting a factor map $\pi:\bfY\to
\sfC^\ast\bfX$ (such exists because by definition $\sfC^\ast\bfX$ is
a factor of some $\sfC$-system).  Now forming $\t{\bfX}:=
\bfX\times_{\{\zeta_{\sfC^\ast}^\bfX = \pi\}}\bfY$, a quick check
using the above fact shows that $\sfC\t{\bfX} = \sfC^\ast\t{\bfX}$,
and that this is equivalent to the $\sfC$-satedness of $\t{\bfX}$.

As remarked previously, a routine argument shows that our basic
examples $\sfZ_0^\L$, $\sfZ_n^\L$ and $\sfZ_{\rm{Ab},n}^\L$ are all
hereditary, and hence that any system is sated with respect to any
of them. However, joins of several hereditary idempotent classes
need not be hereditary.  For example, let $\bfX =
(\bbT^2,\rm{Haar},R_{(\a,0)},R_{(0,\a)})$, where $R_\b$ denotes a
rotation by $\b \in \bbT^2$ and $\a \in (\bbR\setminus \bbQ)/\bbZ$,
and let $\bf{e}_1,\bf{e}_2$ be the standard basis vectors of
$\bbZ^2$. Then clearly $\zeta_0^{T^{\bf{e}_i}} = \pi_{3-i}$, the
projection onto the $(3-i)^{\rm{th}}$ coordinate, and so
$\zeta_0^{T^{\bf{e}_1}}\vee \zeta_0^{T^{\bf{e}_2}} \simeq \id_X$ and
therefore $\bfX$ is itself a system in the class
$\sfZ_0^{\bf{e}_1}\vee \sfZ_0^{\bf{e}_2}$. However,
$\sfZ_0^{\bf{e}_1+\bf{e}_2}\bfX$ is the factor generated by the
SW-NE diagonal circles in $\bbT^2$, so
$\sfZ_0^{\bf{e}_1+\bf{e}_2}\bfX$ is a factor of a system of class
$\sfZ_0^{\bf{e}_1}\vee \sfZ_0^{\bf{e}_2}$ but $\sfZ_0^{\bf{e}_1}\vee
\sfZ_0^{\bf{e}_2}(\sfZ_0^{\bf{e}_1+\bf{e}_2}\bfX)$ is the trivial
system, so $\sfZ_0^{\bf{e}_1+\bf{e}_2}\bfX$ is not
$(\sfZ_0^{\bf{e}_1}\vee \sfZ_0^{\bf{e}_2})$-sated. (Nevertheless,
some preliminary results in Section 7 of~\cite{Aus--ergdirint}
indicate that such counterexamples must be rather special, and we
suspect that the machinery of that paper and the present one may
have more to say on this question in the future.) \fin

It will serve us well to adopt a special name for a particular case
of the above multiple satedness that will recur frequently.

\begin{dfn}[Full isotropy-satedness]\label{dfn:FIS}
A system $T:\bbZ^d \curvearrowright (X,\mu)$ is \textbf{fully
isotropy-sated} (\textbf{FIS}) if whenever $p_i:\bbZ^{r_i}\into
\bbZ^d$, $i=1,2,\ldots,k$, are isomorphic embeddings then the system
$(X,\mu,T)$ is $(\sfZ_0^{p_1}\vee \sfZ_0^{p_2}\vee\cdots\vee
\sfZ_0^{p_k})$-sated.
\end{dfn}

\begin{cor}\label{cor:FIS}
Any $\bbZ^d$-system admits an FIS extension. \qed
\end{cor}

\subsection{Subactions and insensitivity of idempotent classes}

Given an l.c.s.c. group and closed subgroup $\L\leq \G$ there is an
obvious forgetful functor $\G\hbox{-}\Sys \to \L\hbox{-}\Sys$,
$\bfX\mapsto \bfX^{\ \uhr\L}$. Some idempotent classes $\sfC$ in
$\G\hbox{-}\Sys$ actually make sense in both categories, in the
sense that $\bfX \in \sfC$ if and only if $\bfX^{\ \uhr\L}\in
\sfC_0$ for some idempotent class $\sfC_0$ of $\L$-systems. Loosely,
these are the classes that are defined in terms of properties
depending only on the subaction $\bfX^{\ \uhr\L}$, the most obvious
example being $\sfZ_0^\L$. Importantly, in this case forming the
maximal $\sfZ_0^\L$-factor of a $\G$-system $\bfX$ and of its
subaction system $\bfX^{\ \uhr\L}$ give measure-theoretically the
same factor space. This is a simple but important phenomenon that we
will need to appeal to later (although we shall sometimes suppress
the distinction between a class $\sfC$ of $\G$-systems and the
corresponding class $\sfC^{\ \uhr\L}$ of $\L$-systems).

\begin{dfn}[Insensitivity of idempotent classes]
An idempotent class $\sfC$ of $\L$-systems is \textbf{insensitive to
the forgetful functor to $\L$-subactions} if whenever $\bfX$ is a
$\G$-system, the factor map $\zeta_\sfC^{\bfX^{\ \uhr\L}}:X\to\sfC
X$ actually intertwines the whole $\G$-action $T$ with some
$\G$-action on $\sfC X$ (equivalently, if the $\s$-subalgebra
$(\zeta_\sfC^{\bfX^{\ \uhr\L}})^{-1}(\S_{\sfC X})$ is globally
$\G$-invariant).

In this case we can naturally extend the definition of the class
$\sfC$ to the category $\G$-$\Sys$ by letting $\sfC\bfX$ be the
$\G$-action with which $T$ is intertwined by $\zeta_\sfC^{\bfX^{\
\uhr\L}}=: \zeta_\sfC^\bfX$, and we will generally use the same
letter $\sfC$ for the idempotent class in either category.

In general, an idempotent class $\sfC$ of $\G$-systems is
\textbf{insensitive to the forgetful functor to $\L$-subactions} if
it is the extension to $\G$-$\Sys$ of an insensitive idempotent
class of $\L$-systems.
\end{dfn}

This notion of insensitivity will be most important to us in view of
its consequences for satedness.  In order to understand these,
however, we must first introduce an important method for building
system extensions with desirable properties.  For this we restrict
to the setting of countable discrete Abelian groups.

Suppose that $\G$ is a countable discrete Abelian group and $\L\leq
\G$ a subgroup, that $\bfX$ is a $\G$-system and that
$\xi:\bfX'=(X',\mu',S')\to\bfX^{\ \uhr\L}$ is an extension of its
$\L$-subaction system.  It can easily happen that there does not
exist an action $T':\G\curvearrowright(X',\mu')$ such that $(T')^{\
\uhr\L} = S'$. However, if we permit ourselves to pass to further
system extensions we can retrieve this situation, and this will be crucial in cases where our analysis of characteristic factors in the
first place gives information only about a sublattice of $\bbZ^d$,
rather than the whole group.

Here we will introduce a particular construction of such a further
extension (although in general we will use only the abstract fact of
the existence of such an extension). We first need a basic result
and definition from elementary group theory.

\begin{dfn}[Remainder map]
If $\L\leq \G$ are as above and $\O\subseteq \G$ is a fundamental
region for the subgroup $\L$ (that is, $\O$ is a subset containing
exactly one member of every coset in $\G/\L$), then there is a
\textbf{remainder map} $\G\to\O:\g\mapsto R(\g)$ with the property
that $\g - R(\g) \in \L$ for all $\g \in \G$.
\end{dfn}

Now suppose that $\xi:\bfX'\to \bfX^{\ \uhr\L}$ is as above.  We
will give a construction of a further extension of $\bfX'$ based on
a similar idea to that underlying the construction of induced group
representations.

Let $\O \subseteq \G$ be a fundamental region as above with
associated remainder map $R$, chosen so that $e \in \O$ where $e$ is
the identity of $G$. This defines uniquely an `integer part' map $\G
\to\L:\g \mapsto \lfloor\g\rfloor := \g - R(\g)$, and it follows at
once that $R(\g + \g') = R(\g + R(\g')) =R(R(\g)+R(\g'))$.

By Rokhlin's Skew-Product Representation (see, for example, Section
3.3 of Glasner~\cite{Gla03}) the extension $\xi$ of $\L$-systems can
be described by a decomposition of $X$ as $\bigcup_{1\leq n \leq
\infty}A_n$ into $T^{\ \uhr\L}$-invariant sets and, letting
$(Y_n,\nu_n)$ be $\{1,2,\ldots,n\}$ with uniform measure for $1 \leq
n < \infty$ and $(Y_\infty,\nu_\infty)$ be $[0,1)$ with Lebesgue
measure, for each $n$ a Borel cocycle $\Phi_n:\L\times A_n \to
\rm{Aut}(Y_n,\nu_n)$, so that
\[(S')^{\g}(x,y) = (T^\g x,\Phi_n(\g,x)(y))\]
for $(x,y) \in A_n\times Y_n$ for all $\g \in \L$.  We adopt the
simplified notations $\Phi(\g,x):= \Phi_n(\g,x)$ and $(Y_x,\nu_x) :=
(Y_n,\nu_n)$ when $x\in A_n$.  For brevity we will simply write $X'
= X\ltimes Y_\bullet$.

We base our construction of $\t{X}$ on a specification of new,
enlarged fibres above each point $x\in X$. For $x\in X$ let $\t{Y}_x
:= \prod_{\omega \in \O} Y_{T^{-\omega}x}$, and define $\t{X} :=
X\ltimes \t{Y}_\bullet$.  Because $\O$ is countable this product can
be given the structure of a standard Borel space in the obvious way.

We define the $\G$-action $\t{T}$ on $\t{X}$ by its Rokhlin
representation.  Writing a typical element of $Y_\bullet^\O$ as
$(y_\omega)_\omega$, we set
\[\t{T}^\g(x,(y_\omega)_\omega) := \big(T^\g x,\big(\Phi(\lfloor R(\omega - \g) + \g\rfloor,T^{-R(\omega - \g)}x)(y_{R(\omega - \g)})\big)_\omega\big)\]
(noting that this is well-defined: if $(y_\omega)_\omega \in
\prod_{\omega \in \O}Y_{T^{-\omega}x}$, then $y_\omega \in
Y_{T^{-\omega}x}$ for all $\omega \in \O$, and so $y_{R(\omega -
\g)}\in Y_{T^{-R(\omega - \g)}x} = Y_{T^{-\omega}(T^\g x)}$, because
the assignment $x\mapsto Y_x$ is $\L$-invariant).

Using our simple identities for $R$ we can compute that
\begin{eqnarray*}
&&\t{T}^{\g_1}(\t{T}^{\g_2}(x,(y_\omega)_\omega)\\ &&=
\t{T}^{\g_1}\big(T^{\g_2}
x,\big(\Phi(\lfloor R(\omega - \g_2) + \g_2\rfloor,T^{-R(\omega - \g_2)} x)(y_{R(\omega - \g_2)})\big)_\omega\big)\\
&&= \Big(T^{\g_1 + \g_2}
x,\\
&&\quad\quad\big(\Phi(\lfloor R(\omega - \g_1) +
\g_1\rfloor,T^{-R(\omega - \g_1)}T^{\g_2}x)\\
&&\quad\quad\quad\quad\quad\quad\big(\Phi(\lfloor R(R(\omega - \g_1) - \g_2) + \g_2\rfloor,T^{-R(R(\omega - \g_1) - \g_2)} x)(y_{R(R(\omega - \g_1) - \g_2)})\big)\big)_\omega\Big)\\
&&= \Big(T^{\g_1 + \g_2}
x,\\
&&\quad\quad\big(\Phi(\lfloor R(\omega - \g_1) +
\g_1\rfloor,T^{-R(\omega - \g_1)}T^{\g_2}x)\\
&&\quad\quad\quad\quad\quad\quad\big(\Phi(\lfloor R(\omega - (\g_1 +
\g_2)) + \g_2\rfloor,T^{-R(\omega - (\g_1 + \g_2))} x)(y_{R(\omega -
(\g_1 + \g_2))})\big)\big)_\omega\Big),
\end{eqnarray*}
and now we note that
\begin{eqnarray*}
T^{\lfloor R(\omega - (\g_1 + \g_2)) + \g_2\rfloor}T^{-R(\omega -
(\g_1 + \g_2))}x &=& T^{\lfloor \omega - (\g_1 + \g_2) + \g_2\rfloor
- \lfloor \omega - (\g_1 + \g_2)\rfloor -R(\omega - (\g_1 +
\g_2))}x\\
&=& T^{\lfloor \omega - \g_1\rfloor - (\omega - (\g_1 + \g_2))}x\\
&=&
T^{(\omega - \g_1) - R(\omega - \g_1) - (\omega - (\g_1 + \g_2))}x\\
&=& T^{-R(\omega - \g_1)}T^{\g_2}x,
\end{eqnarray*}
and so the cocycle equation for $\Phi$ gives
\begin{eqnarray*}
&&\Phi(\lfloor R(\omega - \g_1) + \g_1\rfloor,T^{-R(\omega -
\g_1)}T^{\g_2}x)\circ \Phi(\lfloor R(\omega - (\g_1 + \g_2)) +
\g_2\rfloor,T^{-R(\omega - (\g_1 + \g_2))} x)\\
&&= \Phi(\lfloor R(\omega - \g_1) + \g_1\rfloor + \lfloor R(\omega -
(\g_1 + \g_2)) + \g_2\rfloor,T^{-R(\omega - (\g_1 + \g_2))} x)\\
&&= \Phi(-\lfloor \omega - \g_1\rfloor + \lfloor \omega -
\g_1\rfloor - \lfloor \omega - (\g_1 + \g_2)\rfloor,T^{-R(\omega -
(\g_1 + \g_2))} x).
\end{eqnarray*}
Inserting this into the above formula for
$\t{T}^{\g_1}(\t{T}^{\g_2}(x,(y_\omega)_\omega))$ shows that it is
equal to $\t{T}^{\g_1 + \g_2}(x,(y_\omega)_{\omega})$, and hence
that $\t{T}$ is a $\G$-action.

Finally, if we let
\[\pi:(x,(y_\omega)_\omega)\mapsto x\quad\quad\hbox{and}\quad\quad \a:(x,(y_\omega)_\omega)\mapsto (x,y_e)\]
(recalling that $e\in \O$) then it is routine to check that
\[\pi(\t{T}^\g(x,(y_\omega)_\omega)) = T^\g x = T^\g \pi(x,(y_\omega)_\omega)\quad\quad \forall \g \in \G\]
and
\begin{multline*}
\a(\t{T}^\g(x,(y_\omega)_\omega)) = \a\big(T^\g x,\big(\Phi(\lfloor
R(\omega - \g) + \g\rfloor,T^{-R(\omega - \g)}x)(y_{R(\omega -
\g)})\big)_\omega\big)\\
= (T^\g,\Phi(\lfloor R( - \g) + \g\rfloor,T^{-R(- \g)}x)(y_{R(-
\g)})) = (T^\g,\Phi(\g,x)(y_e)) = (S')^\g(x,y_e)
\end{multline*}
when $\g \in \L$, so we have the required commutative diagram
\begin{center}
$\phantom{i}$\xymatrix{ \t{\bfX}^{\ \uhr\L}\ar[rr]^\pi\ar[dr]^\a & & \bfX^{\ \uhr\L}\\
& \bfX'\ar[ur]^\xi. }
\end{center}

\begin{dfn}[Fibrewise power extension]\label{dfn:FP}
We will refer to the particular extension $\pi:\t{\bfX}\to\bfX$
constructed above as the \textbf{fibrewise power extension} (or
\textbf{FP extension}) of $\bfX$ corresponding to the subgroup
$\L\leq \G$ and system extension $\xi$.
\end{dfn}

\textbf{Remark}\quad Interestingly, the appeal made to the
discreteness of $\G/\L$ in the above proof seems to be quite
important.  While other instances of this theorem are certainly
available, it seems to be difficult to prove a comparably general
statement for an inclusion $\L \leq \G$ of arbitrary locally
compact second countable Abelian groups; it would be interesting to
know whether some alternative construction could be found to handle
that setting. It is also worth remarking that there are certainly
pairs of locally compact non-Abelian groups for which the conclusion
fails: for example, by the Howe-Moore Theorem~(see, for instance,
Section 3.3 of \cite{FerKat02}) any ergodic action of a non-compact
connected simple Lie group $G$ with finite centre is mixing, and so
any non-mixing ergodic action of $\bbR$ cannot be extended to a
larger $\bbR$-system in which the action can be enlarged to the
whole group $G$ for any embedding $\bbR \into G$ as a one-parameter
subgroup. \fin

\textbf{Remark}\quad Since the submission of the present paper, the above ideas have been simplified and generalized to cover any inclusion of countable acting groups: see~\cite[Theorem 2.1]{Aus--commute}. \fin

We can now quickly derive some consequences for satedness.

\begin{lem}\label{lem:sated-insens}
If $\L\leq \G$ are as above and $\sfC$ is an idempotent class that
is insensitive to the forgetful functor to $\L$-actions then a
$\G$-system $\bfX$ is $\sfC$-sated if and only if $\bfX^{\ \uhr\L}$
is $\sfC$-sated.
\end{lem}

\textbf{Proof}\quad It is clear that if $\bfX$ admits a
$\sfC$-adjoining that is not relatively independent over $\sfC\bfX$,
then simply applying the forgetful functor gives the same phenomenon
among the $\L$-subactions: thus the $\sfC$-satedness of $\bfX^{\
\uhr\L}$ implies that of $\bfX$.  The reverse direction follows
similarly, except that a given $\sfC$-adjoining of $\bfX^{\ \uhr\L}$
may need to be extended further (for example, to an FP extension) to
recover an action of the whole of $\G$, and this extension of $\bfX$
will then witness that it is not sated. \qed

\begin{cor}\label{cor:isotropy-insens}
For any subgroups $\L_1,\L_2,\ldots,\L_k \leq \bbZ^d$ a system
$\bfX$ is $\big(\bigvee_{i\leq k}\sfZ_0^{\L_i}\big)$-sated if and
only if $\bfX^{\ \uhr(\L_1 + \L_2 + \ldots + \L_k)}$ is
$\big(\bigvee_{i\leq k}\sfZ_0^{\L_i}\big)$-sated. \qed
\end{cor}

\section{Some applications to characteristic factors}

\subsection{The Furstenberg self-joining}

Consider a $\bbZ^d$-system $\bfX = (X,\mu,T)$. As have many previous
works in this area, our analysis of characteristic factors
associated to different sequences of linear averages will make heavy
use of a particular class of self-joinings of $\bfX$. Given $k$
directions $\bf{p}_1$, $\bf{p}_2$, \ldots, $\bf{p}_k\in \bbZ^d$, let
us here write
\[S_N(f_1,f_2,\ldots,f_k):= \frac{1}{N}\sum_{n=1}^N\prod_{i=1}^kf_i\circ T^{n\bf{p}_i}\]
for the associated $k$-fold linear nonconventional averages. Now for
$A_1$, $A_2$, \ldots, $A_k \in \S_X$ we can define
\[\mu^{\rm{F}}_{T^{\bf{p}_1},T^{\bf{p}_2},\ldots,T^{\bf{p}_k}}(A_1\times A_2\times\cdots\times A_k) := \lim_{N\to\infty}\int_XS_N(1_{A_1},1_{A_2},\ldots,1_{A_k})\,\d\mu,\]
where the existence of this limit follows from the known
convergence of linear nonconventional averages. (Once convergence of
the relevant polynomial nonconventional averages has been
established, a similar definition can be made corresponding to such
polynomial averages, but in the nonlinear case these have yet to
prove similarly useful.)

Now it is routine to show (recalling that $X$ is standard Borel)
that the above definition extends by multilinearity and continuity to a $k$-fold
self-joining of $\mu$ on $X^k$, which is invariant under not only
the $\bbZ^d$-action $T^{\times k}$ but also the `diagonal
transformation' $\vec{T} := T^{\bf{p}_1}\times T^{\bf{p}_2}\times
\cdots\times T^{\bf{p}_k}$.  This is the \textbf{Furstenberg
self-joining} of $\mu$ associated to the transformations
$T^{\bf{p}_1}$, $T^{\bf{p}_2}$,\ldots, $T^{\bf{p}_k}$, and will be
denoted by
$\mu^{\rm{F}}_{T^{\bf{p}_1},T^{\bf{p}_2},\ldots,T^{\bf{p}_k}}$ or
$\mu^{\rm{F}}_{\vec{T}}$, or sometimes abbreviated to
$\mu^{\rm{F}}$. Clearly whenever $f_i \in L^\infty(\mu)$ for
$i=1,2,\ldots,k$ we have also
\[\int_{X^k}f_1\otimes f_2\otimes\cdots\otimes
f_k\,\d\mu^{\rm{F}}_{T^{\bf{p}_1},T^{\bf{p}_2},\ldots,T^{\bf{p}_k}}
= \lim_{N\to\infty}\int_X S_N(f_1,f_2,\ldots,f_k)\,\d\mu.\]

In fact, in~\cite{Aus--nonconv} the convergence of the averages that define the Furstenberg self-joining is proved alongside the convergence of the functional nonconventional averages themselves as part of a zigzag induction (one claim for a given $k$ implies the
other for that $k$, which then implies the first for $k+1$, and so
on). This is possible in view of a reduction of the above limits to
the study of linear averages involving only $k-1$ transformations,
which will also be important for us here: simply because $\mu$ is
$T$-invariant, and working now in terms of bounded functions rather
than sets, we can re-write the above limit as
\begin{multline*}
\int_{X^k}f_1\otimes f_2\otimes\cdots\otimes
f_k\,\d\mu^{\rm{F}}_{T^{\bf{p}_1},T^{\bf{p}_2},\ldots,T^{\bf{p}_k}}\\
= \int_Xf_1\cdot\Big(\lim_{N\to\infty}\frac{1}{N}\sum_{n=1}^N
\prod_{i=2}^k(f_i\circ T^{n(\bf{p}_i - \bf{p}_1)})\Big)\,\d\mu.
\end{multline*}

Knowing that the Furstenberg self-joining
$\mu^{\rm{F}}_{T^{\bf{p}_1},T^{\bf{p}_2},\ldots,T^{\bf{p}_k}}$
exists, a basic application of the Hilbert space version of the
classical van der Corput estimate gives us a way to use it to
control the asymptotic behaviour of $S_N(f_1,f_2,\ldots,f_k)$, in
the sense of the following estimate taken from Lemma 4.7
of~\cite{Aus--nonconv}. (We have modified the statement to give
explicit bounds, but the proof is unchanged.)

\begin{lem}\label{lem:Fberg-controls-aves}
If $f_1$, $f_2$, \ldots , $f_k \in L^\infty(\mu)$ are $1$-bounded, then there is a $1$-bounded $\vec{T}$-invariant function $g\in L^\infty(\mu^{\rm{F}}_{\vec{T}})$ such that
\[\Big|\int_{X^k}\prod_{i=1}^k(f_i\circ\pi_i)\cdot g\,\d\mu^{\rm{F}}_{\vec{T}}\Big| \geq \lim_{N\to\infty}\|S_N(f_1,f_2,\ldots,f_d)\|^2_2.\]
\qed
\end{lem}

This now implies a useful sufficient condition for characteristicity
of a tuple of factors.

\begin{cor}\label{cor:Fberg-controls-charfactors}
A tuple of factors $\xi_i:\bfX\to\bfY_i$, $i=1,2,\ldots,k$, is
characteristic for the averages $S_N$ if for any choice of $f_i\in
L^\infty(\mu)$, $i=1,2,\ldots,k$, and $\vec{T}$-invariant $g\in
L^\infty(\mu^{\rm{F}}_{\vec{T}})$ we have
\[\int_{X^k}\prod_{i=1}^k(f_i\circ\pi_i)\cdot g\,\d\mu^{\rm{F}}_{\vec{T}} = \int_{X^k}\prod_{i=1}^k(\sfE_\mu(f_i\,|\,\xi_i)\circ\pi_i)\cdot g\,\d\mu^{\rm{F}}_{\vec{T}}.\]
\end{cor}

\textbf{Proof}\quad Suppose that $f_i \in L^\infty(\mu)$ for each
$i=1,2,\ldots,k$.  We need to show that
\[S_N(f_1,f_2,\ldots,f_k)\sim S_N(\sfE_\mu(f_1\,|\,\xi_1),\sfE_\mu(f_2\,|\,\xi_2),\ldots,\sfE_\mu(f_k\,|\,\xi_k))\]
as $N\to\infty$, but by replacing each function with its conditional
expectation in turn it clearly suffices to show that
\[S_N(f_1,f_2,\ldots,f_k)\sim S_N(\sfE_\mu(f_1\,|\,\xi_1),f_2,\ldots,f_k).\]
This, in turn, is equivalent to
\[S_N(f_1 - \sfE_\mu(f_1\,|\,\xi_1),f_2,\ldots,f_k)\to 0,\]
and this now follows from the assumption and
Lemma~\ref{lem:Fberg-controls-aves} because for any
$\vec{T}$-invariant $g\in L^\infty(\mu^{\rm{F}}_{\vec{T}})$ we have
\begin{multline*}
\int_{X^k}\prod_{i=1}^k(f_i\circ\pi_i)\cdot
g\,\d\mu^{\rm{F}}_{\vec{T}} =
\int_{X^k}\prod_{i=1}^k(\sfE_\mu(f_i\,|\,\xi_i)\circ\pi_i)\cdot
g\,\d\mu^{\rm{F}}_{\vec{T}}\\
=\int_{X^k}\sfE_\mu(f_1\,|\,\xi_i)\cdot\prod_{i=2}^k(f_i\circ\pi_i)\cdot
g\,\d\mu^{\rm{F}}_{\vec{T}}
\end{multline*}
and so
\[\int_{X^d}(f_1 - \sfE_\mu(f_1\,|\,\xi_1))\cdot\prod_{i=2}^d(f_i\circ\pi_i)\cdot g\,\d\mu^{\rm{F}}_{\vec{T}} = 0.\]
\qed

\textbf{Remark}\quad An alternative to the Furstenberg self-joining
that can sometimes be put to similar uses has been constructed by
Host and Kra, first in the case of powers of a single transformation
in~\cite{HosKra05} and then for several commuting transformations by
Host in~\cite{Hos09}.  This is defined in terms of a tower of
iterated relatively independent self-products, and so has the
advantage over the Furstenberg self-joining that it does not require
an appeal to a previously-known nonconventional convergence result
for its definition.  Although we shall focus on the Furstenberg self-joining here for consistency, I suspect that the present paper and its sequel could be re-worked to use a Host-Kra self-joining throughout, and that neither presentation would be substantially easier. (In early drafts of these papers, the preference for the Furstenberg self-joining was dictated by a particular appeal to it in the last stages of proving Theorem~\ref{thm:polyconv}, but subsequent improvements to that proof have made these unnecessary.) \fin

It is easy to see that a tuple of factors $(\xi_1, \xi_2,
\ldots, \xi_d)$ is characteristic for the averages $S_N$ if and only
if each of the $d$ tuples
\begin{eqnarray*}
&&(\xi_1,\id_X,\id_X,\ldots,\id_X),\\
&&(\id_X,\xi_2,\id_X,\ldots,\id_X),\\
&&\quad\quad\quad\quad\quad\quad\vdots\\
&&(\id_X,\id_X,\id_X,\ldots,\xi_d)
\end{eqnarray*}
is characteristic for them.  A slightly more subtle property that we
will find useful later is the following.

\begin{lem}\label{lem:char-factors-symmetry}
For any factor $\xi:\bfX\to \bfY$ the tuple
$(\xi,\id_X,\ldots,\id_X)$ is characteristic for the nonconventional
averages
\[\frac{1}{N}\sum_{n=1}^N\prod_{i=1}^kf_i\circ T^{n\bf{p}_i}\quad\quad f_1,f_2,\ldots,f_k\in L^\infty(\mu)\]
if and only if the tuple $(\id_X,\xi,\id_X,\ldots,\id_X)$ is
characteristic for the nonconventional averages
\[\frac{1}{N}\sum_{n=1}^N(f_0\circ T^{-n\bf{p}_j})\cdot\prod_{i\leq k,\,i\neq j}f_i\circ T^{n(\bf{p}_i - \bf{p}_j)}\quad\quad f_0,f_1,\ldots,f_{j-1},f_{j+1},\ldots,f_k \in L^\infty(\mu)\]
for every $j = 2,3,\ldots,k$.
\end{lem}

\textbf{Proof}\quad This follows from a similar re-arrangement to
those we have already seen above.  By symmetry it suffices to treat
only one of the needed implications, so let us suppose that
$(\xi,\id_X,\ldots,\id_X)$ is characteristic for
\[\frac{1}{N}\sum_{n=1}^N\prod_{i=1}^kf_i\circ T^{n\bf{p}_i}\]
and show that $(\id_X,\xi,\id_X,\ldots,\id_X)$ is characteristic for
\[\frac{1}{N}\sum_{n=1}^N(f_0\circ T^{-n\bf{p}_2})\cdot\prod_{i\leq k,\,i\neq 2}f_i\circ T^{n(\bf{p}_i - \bf{p}_2)}\quad\quad f_0,f_1,f_3,\ldots,f_k \in L^\infty(\mu).\]
Replacing $f_1$ by $f_1 - \sfE_\mu(f_1\,|\,\xi_1)$, it will suffice
to show that if the latter averages do not tend to zero in
$L^2(\mu)$ for some choice of $f_0,f_3,\ldots,f_k$ then also the
former do not tend to zero for some choice of $f_2,f_3,\ldots,f_k$.
Thus, suppose that $f_0,f_3,\ldots,f_k$ are such that the latter
averages do not tend to zero, and now let
\[f_2 := \lim_{N\to\infty}\frac{1}{N}\sum_{n=1}^N(f_0\circ T^{-n\bf{p}_2})\cdot\prod_{i\leq k,\,i\neq 2}f_i\circ T^{n(\bf{p}_i - \bf{p}_2)}.\]
The condition that $f_2 \neq 0$ and a change of variables now give
\begin{eqnarray*}
0&\neq&\lim_{N\to\infty}\int_Xf_2\cdot
\Big(\frac{1}{N}\sum_{n=1}^N(f_0\circ
T^{-n\bf{p}_2})\cdot\prod_{i\leq k,\,i\neq 2}f_i\circ T^{n(\bf{p}_i
- \bf{p}_2)}\Big)\,\d\mu\\
&&= \lim_{N\to\infty}\frac{1}{N}\sum_{n=1}^N\int_Xf_2\cdot(f_0\circ
T^{-n\bf{p}_2})\cdot\prod_{i\leq k,\,i\neq 2}f_i\circ T^{n(\bf{p}_i
- \bf{p}_2)}\,\d\mu\\
&&= \lim_{N\to\infty}\frac{1}{N}\sum_{n=1}^N\int_Xf_0\cdot(f_2\circ
T^{n\bf{p}_2})\cdot\prod_{i\leq k,\,i\neq 2}f_i\circ
T^{n\bf{p}_i}\,\d\mu\\
&&= \lim_{N\to\infty}\int_X
f_0\cdot\Big(\frac{1}{N}\sum_{n=1}^N\prod_{i=1}^kf_i\circ
T^{n\bf{p}_i}\Big)\,\d\mu,
\end{eqnarray*}
and so we must also have
\[\frac{1}{N}\sum_{n=1}^N\prod_{i=1}^kf_i\circ T^{n\bf{p}_i}\not\to 0,\]
as required. \qed

\textbf{Example}\quad An easily-generalized argument of
Ziegler~\cite{Zie07} in the case $d=1$ and $\bf{p}_i = a_i \in \bbZ$
shows that there is always a unique minimal characteristic factor
tuple: a characteristic tuple of factors $\xi_i:\bfX\to\bfY_i$,
$i=1,2,\ldots,k$, such that any other characteristic tuple of
factors $\xi'_i:\bfX\to\bfY'_i$ must satisfy $\xi'_i \succsim \xi_i$
for all $i\leq k$.  However, it is worth noting that the members of
this tuple can depend on the whole system $\bfX$, in that if we
restrict to averages involving functions $f_j$ that are all lifted
from $\bfY_i$ for some fixed $i$, then $\bfY_i$ may in turn admit a
characteristic tuple of smaller factors.

For example, when $d=2$ and $\bf{p}_i = \bf{e}_i$ for $i=1,2$,
consider three irrational and rationally independent points on the
circle $r,s,t \in \bbT$, and let $\bfX$ be the $\bbZ^2$ system on
$(X,\mu) = (\bbT^2,m_{\bbT^2})$ generated by $T_1 := R_s\times R_t$
and $T_2 := R_r\times R_s$.  In this simple setting we can use
Fourier analysis to obtain that the minimal characteristic factors
$\xi_1$, $\xi_2$ are equivalent to the first and second coordinate
projections $\bbT^2\to\bbT$ respectively. However, after passing
down through the first coordinate projection, it is equally easy to
compute that the minimal characteristic factors of the resulting
system on $(\bbT,m_\bbT)$ are both trivial.

It follows that there is in general no tuple of idempotent classes
of systems $\sfC_1$, $\sfC_2$, \ldots, $\sfC_k$ such that the factors $\bfY_i =
\sfC_i\bfX$ serve as the minimal characteristic tuple of factors for every $\bfX$.  This contrasts interestingly with the
case $d=1$, where the main technical result of Host and
Kra~\cite{HosKra05} and Ziegler~\cite{Zie07} can be phrased as
asserting that there is such a class, and for $k$ distinct integers
$p_1$, $p_2$, \ldots, $p_k\in \bbZ$ we have $\sfC_1 = \sfC_2 =
\ldots = \sfC_k$ and it is the class of all `direct integrals'
(suitably defined) of inverse limits of $k$-step nilsystems. \fin

Let us finish by recording the following useful property of minimal
characteristic factors.

\begin{lem}\label{lem:char-factor-inv-lim}
If $(\bfX_{(m)})_{m\geq 0}$, $(\psi^{(m)}_{(k)})_{m\geq k\geq 0}$ is
an inverse system with inverse limit $\bfX_{(\infty)}$,
$(\psi_{(m)})_{m\geq 0}$ and the factors of the minimal
characteristic tuples of these systems for some averaging scheme are
$\xi_{(m),i}$, $i=1,2,\ldots,k$ and $\xi_{(\infty),i}$,
$i=1,2,\ldots,k$ respectively then
\[\xi_{(\infty),i} = \bigvee_{m\geq 0}\xi_{(m),i}\circ\psi_{(m)}.\]
\end{lem}

\textbf{Proof}\quad The direction $\succsim$ is obvious (since any
particular nonconventional averages on system $\bfX_{(m)}$ can be
lifted to $\bfX_{(\infty)}$), so we need only show the reverse
containment.

To this end, suppose that $f_i \in L^\infty(\mu_{(\infty)})$ for
$i=1,2,\ldots,k$.  Then by the definition of the inverse limit, we
know that we can approximate these functions arbitrarily well in
$L^2(\mu)$ by functions of the form $g_i\circ\psi_{(m)}$ for $g_i
\in L^\infty(\mu_{(m)})$ with $\|g_i\|_\infty \leq \|f_i\|_\infty$
and $m$ sufficiently large. This approximation now clearly gives
\[S_{(\infty),N}(f_1,f_2,\ldots,f_k) \approx S_{(m),N}(g_1,g_2,\ldots,g_k)\circ\psi_{(m)}\quad\quad\hbox{in }L^2(\mu)\]
uniformly in $N$, and this latter behaves asymptotically as
\begin{eqnarray*}
&&S_{(m),N}(\sfE_{\mu_{(m)}}(g_1\,|\,\xi_{(m),1}),\sfE_{\mu_{(m)}}(g_2\,|\,\xi_{(m),2}),\ldots,\sfE_{\mu_{(m)}}(g_k\,|\,\xi_{(m),k}))\circ\psi_{(m)}\\
&&=
S_{(\infty),N}(\sfE_{\mu_{(\infty)}}(g_1\circ\psi_{(m)}\,|\,\xi_{(m),1}\circ\psi_{(m)}),\sfE_{\mu_{(\infty)}}(g_2\circ\psi_{(m)}\,|\,\xi_{(m),2}\circ\psi_{(m)}),\\
&&\quad\quad\quad\quad\quad\quad\quad\quad\quad\quad\quad\quad\quad\quad\quad\quad
\ldots,\sfE_{\mu_{(\infty)}}(g_k\circ\psi_{(m)}\,|\,\xi_{(m),k}\circ\psi_{(m)}))
\end{eqnarray*}
in $L^2(\mu)$ as $N\to\infty$, by the defining property of
$\xi_{(m),1}$, $\xi_{(m),2}$, \ldots, $\xi_{(m),k}$.

This, in turn, is approximately equal to
\[S_{(\infty),N}(\sfE_{\mu_{(\infty)}}(f_1\,|\,\xi_{(m),1}\circ\psi_{(m)}),\sfE_{\mu_{(\infty)}}(f_2\,|\,\xi_{(m),2}\circ\psi_{(m)}),\ldots,\sfE_{\mu_{(\infty)}}(f_k\,|\,\xi_{(m),k}\circ\psi_{(m)})),\]
and so letting $m\to\infty$ and observing that each
$(\sfE_{\mu_{(\infty)}}(f_i\,|\,\xi_{(m),i}\circ\psi_{(m)}))_{m\geq
1}$ for $i=1,2,\ldots,d$ is a uniformly bounded martingale, we
obtain
\begin{multline*}
S_{(\infty),N}(f_1,f_2,\ldots,f_k)\\
\sim
S_{(\infty),N}(\sfE_{\mu_{(\infty)}}(f_1\,|\,\xi^\circ_1),\sfE_{\mu_{(\infty)}}(f_2\,|\,\xi^\circ_2),\ldots,\sfE_{\mu_{(\infty)}}(f_k\,|\,\xi^\circ_k))
\end{multline*}
as $N\to\infty$ with
\[\xi^\circ_i := \bigvee_{m\geq 0}\xi_{(m),i}\circ\psi_{(m)},\]
and hence $\xi_{(\infty),i}\simeq \xi^\circ_i$, as required. \qed

\subsection{Linearly independent directions in discrete time}\label{subs:roots}

In this section we will address the easier of the questions posed in
the introduction: whether we can construct pleasant extensions while
retaining the existence of roots for our transformations.  In fact
it will follow quite easily from the machinery developed above that
FIS extensions achieve this goal (Definition~\ref{dfn:FIS}).

\begin{prop}[FIS extensions are pleasant]\label{prop:FIS-implies-pleasant}
If $(X,\mu,T)$ is an FIS system and $\bf{p}_1$, $\bf{p}_2$, \ldots,
$\bf{p}_k \in \bbZ^d$ are linearly independent then the tuple of
factors
\[\xi_i := \zeta_0^{T^{\bf{p}_i}}\vee\bigvee_{j\in \{1,2,\ldots,k\}\setminus \{i\}}\zeta_0^{T^{\bf{p}_i} = T^{\bf{p}_j}}\quad\quad i=1,2,\ldots,k\] is characteristic
for the associated linear nonconventional averages
\[\frac{1}{N}\sum_{n=1}^N\prod_{i=1}^kf_i\circ T^{n\bf{p}_i}.\]
\end{prop}

The case in which in addition $\bbZ\bf{p}_1 + \bbZ\bf{p}_2 + \cdots
+ \bbZ\bf{p}_k = \bbZ^d$ (so necessarily $k = d$) is implicitly
contained in~\cite{Aus--nonconv}; the point here is to handle the
case when the $\bf{p}_j$ generate a proper sublattice.

Before turning to Proposition~\ref{prop:FIS-implies-pleasant} we
illustrate our basic method by proving the following useful lemma
(which is, in turn, implicitly contained
in~\cite{Aus--newmultiSzem}).  For our purposes a tuple of
isomorphic embeddings $p_j:\bbZ^{r_j}\into \bbZ^d$, $j=1,2,\ldots,k$
is \textbf{totally linearly independent} if
\[p_1(\bf{n}_1) + p_2(\bf{n}_2) + \ldots + p_k(\bf{n}_k) = \bs{0}\in\bbZ^d\quad\quad\Rightarrow\quad\quad \bf{n}_j = \bs{0}\in \bbZ^{r_j}\ \forall j\leq k.\]

\begin{lem}
If $(X,\mu,T)$ is an FIS system then whenever
$p_j:\bbZ^{r_j}\into\bbZ^d$ are totally linearly independent
isomorphic embeddings for $j=1,2,\ldots,k$ we have
\[\zeta_0^{T^{p_1}}\wedge \Big(\bigvee_{j=2}^k\zeta_0^{T^{p_j}}\Big) \simeq \bigvee_{j=2}^k\zeta_0^{T^{p_1\oplus p_j}},\]
where $T^{p_j}$ is the $\bbZ^{r_j}$-action $\bf{n}\mapsto
T^{p_j(\bf{n})}$.
\end{lem}

\textbf{Proof}\quad Let $\L :=
p_1(\bbZ^{r_1})+p_2(\bbZ^{r_2})+\cdots +p_k(\bbZ^{r_k}) \leq
\bbZ^d$. We first suppose $\L = \bbZ^d$, and then use this to prove
the general case. Note that since our $p_j$ are injective and
linearly independent, in this special case they together define an
isomorphism $\bbZ^d \cong \bbZ^{r_1}\oplus \bbZ^{r_2}\oplus
\cdots\oplus \bbZ^{r_k}$.

It is clear that for any system $(X,\mu,T)$ we have
\[\zeta_0^{T^{p_1}}\wedge \Big(\bigvee_{j=2}^k\zeta_0^{T^{p_j}}\Big) \succsim \bigvee_{j=2}^k\zeta_0^{T^{p_1\oplus p_j}},\]
so we need only prove the reverse containment. Let
\[(\t{X},\t{\mu})
:=
(X\times_{\zeta_0^{T^{p_1}}}X,\mu\otimes_{\zeta_0^{T^{p_1}}}\mu)\]
and $\pi_1$ and $\pi_2$ be respectively the first and second
coordinate projections $\t{X} \to X$, and define a $\bbZ^d$-action
$\t{T}$ on $(\t{X},\t{\mu})$ by setting
\[\t{T}^{p_i(\,\cdot\,)} := \left\{\begin{array}{ll}T^{p_1(\,\cdot\,)}\times \id_X&\quad\quad\hbox{if }i=1\\ (T^{\times 2})^{p_i(\,\cdot\,)}&\quad\quad\hbox{if }i=2,3,\ldots,k.\end{array}\right.\]
and extending additively. This is easily to seen to be a
well-defined probability-preserving $\bbZ^d$-system and an extension
of $(X,\mu,T)$ through $\pi_1$. Now note that the whole second
coordinate in $\t{X}$ is $\t{T}^{p_1}$-invariant, and hence that
\begin{multline*}
\bigvee_{j=2}^k\zeta_0^{\t{T}^{p_1\oplus p_j}}\succsim
\Big(\bigvee_{j=2}^k\zeta_0^{T^{p_j}}\Big)\circ \pi_2\\ \succsim
\Big(\zeta_0^{T^{p_1}}\wedge\bigvee_{j=2}^k\zeta_0^{T^{p_j}}\Big)\circ\pi_2
\simeq
\Big(\zeta_0^{T^{p_1}}\wedge\bigvee_{j=2}^k\zeta_0^{T^{p_j}}\Big)\circ\pi_1,
\end{multline*}
where the last equivalence holds because
$\zeta_0^{T^{p_1}}\circ\pi_1 \simeq \zeta_0^{T^{p_1}}\circ\pi_2$ by
construction. On the other hand since $(X,\mu,T)$ is FIS the factors
$\bigvee_{j=2}^k\zeta_0^{\t{T}^{p_1\oplus p_j}}$ and $\pi_1$ must be
relatively independent over
\[\pi_1\wedge
\Big(\bigvee_{j=2}^k\zeta_0^{\t{T}^{p_1\oplus p_j}}\Big) \simeq
\Big(\bigvee_{j=2}^k\zeta_0^{T^{p_1\oplus p_j}}\Big)\circ\pi_1,\] so
in fact we have
\[\Big(\bigvee_{j=2}^k\zeta_0^{T^{p_1\oplus p_j}}\Big)\circ\pi_1\succsim \Big(\zeta_0^{T^{p_1}}\wedge\bigvee_{j=2}^k\zeta_0^{T^{p_j}}\Big)\circ\pi_1\]
and hence
\[\bigvee_{j=2}^k\zeta_0^{T^{p_1\oplus p_j}}\succsim \zeta_0^{T^{p_1}}\wedge\bigvee_{j=2}^k\zeta_0^{T^{p_j}},\] as required.

Finally, for a general $\L$ Corollary~\ref{cor:isotropy-insens}
tells us that the subaction system $\bfX^{\ \uhr\L}$ is still FIS,
and so since all joins of the idempotent classes $\sfZ_0^{p_i}$ are
insensitive to the forgetful functor to $\L$-subactions the special
case treated above completes the proof. \qed

\textbf{Proof of Proposition~\ref{prop:FIS-implies-pleasant}}\quad
Let $\L := \bbZ\bf{p}_1 + \bbZ\bf{p}_2 + \cdots + \bbZ\bf{p}_k$.

Once again we first treat the case $\L = \bbZ^d$; this is already
covered in~\cite{Aus--nonconv} in slightly different terms, and the
underlying idea here is as in that paper. Write $T_i :=
T^{\bf{p}_i}$ and fix some $j \leq k$.  Consider the extension
$\pi:\t{\bfX}\to\bfX$ built from the Furstenberg self-joining by
\begin{itemize}
\item letting $(\t{X},\t{\mu}) := (X^k,\mu^{\rm{F}})$,
\item defining the lifted transformations $\t{T}_i$ by
\[\t{T}_i = \left\{\begin{array}{ll} T_1\times T_2\times \cdots\times T_k&\quad\quad\hbox{for }i=j\\ T_i^{\times k}&\quad\quad\hbox{for }i\in\{1,2,3,\ldots,k\}\setminus\{j\},\end{array}\right.\]
\item writing $\pi_i:X^k\to X$, $i=1,2,\ldots,k$, for the coordinate
projections,
\item and taking $\pi := \pi_j$.
\end{itemize}

Now let $f_i \in L^\infty(\mu)$ for $i=1,2,\ldots,k$ and let $g \in
L^\infty(\mu^{\rm{F}})$ be $\t{T}_j$-invariant. Observe from the
above choice of the lifted transformations that $f_i\circ\pi_i$ is a
$\t{T}_i\t{T}_j^{-1}$-invariant function on $\t{X}$ for each $i\neq
j$, and so the function $\prod_{i\leq k,\,i\neq
j}(f_i\circ\pi_i)\cdot g$ on $\t{X}$ is
$\big(\zeta_0^{\t{T}_j}\vee\bigvee_{i\leq k,\,i\neq
j}\zeta_0^{\t{T}_i = \t{T}_j}\big)$-measurable.  Since $\bfX$ is
FIS, under $\mu^{\rm{F}}$ this function is relatively independent
from $f_j\circ\pi_j$ over $\big(\zeta_0^{T_j}\vee\bigvee_{i\leq
k,\,i\neq j}\zeta_0^{T_i = T_j}\big)\circ\pi_j$: that is, writing
$\xi_j := \zeta_0^{T_j}\vee\bigvee_{i\leq k,\,i\neq j}\zeta_0^{T_i =
T_j}$, we have
\[\int_{X^k} \prod_{i=1}^k(f_i\circ\pi_i)\cdot g\,\d\mu^{\rm{F}} = \int_{X^k} (\sfE_\mu(f_j\,|\,\xi_j)\circ\pi_j)\cdot\prod_{i\leq k,\,i\neq j}^k(f_i\circ\pi_i)\cdot g\,\d\mu^{\rm{F}}.\]
Using this argument to replace $f_j$ with $\sfE_\mu(f_j\,|\,\xi_j)$
for each $j$ in turn, Corollary~\ref{cor:Fberg-controls-charfactors}
tells us that $(\xi_1, \xi_2,\ldots,\xi_k)$ is characteristic, as
required.

Now if $\L$ is a general sublattice, we observe that for a given
tuple of factors of our FIS system $\bfX$, their characteristicity
depends only on the subaction system $\bfX^{\ \uhr\L}$, which is
also FIS by Corollary~\ref{cor:isotropy-insens}, and so as for the
preceding lemma the above special case completes the proof. \qed

\subsection{An example in continuous time}

In addition to the above description of pleasant extensions for
certain linear averages (by itself only a very modest generalization
of technical results from~\cite{Aus--nonconv}), we will now offer an
application of sated extensions to a different convergence problem
for nonconventional averages.  This problem is `quadratic' and
`two-dimensional', which features introduce new difficulties, but it
is also in `continuous time', and we will find that this allows us
to recover a fairly short proof.

Thus, we now switch to the setting of a jointly measurable action
$\bbR^2\curvearrowright (X,\mu)$, which we denote by $\bbR^2\to
\rm{Aut}_0(X,\mu):\bf{v} \mapsto \tau^{\bf{v}}$.  We also let
$\bf{e}_1$, $\bf{e}_2$ be the standard basis of $\bbR^2$.

\begin{thm}\label{thm:cts-poly-aves}
The averages
\[S_T(f_1,f_2):= \frac{1}{T}\int_0^T (f_1\circ \tau^{t^2\bf{e}_1})(f_2\circ \tau^{t^2\bf{e}_1 + t\bf{e}_2})\,\d t\]
converge in $L^2(\mu)$ as $T \to\infty$ for any $f_1,f_2 \in
L^\infty(\mu)$.
\end{thm}

As in~\cite{Aus--nonconv} this will follow once we ascend to a
suitable pleasant extension.

\begin{prop}[Pleasant extensions for continuous-time quadratic averages]\label{prop:char-cts-poly}
If the $\bbR^2$-system $(X,\mu,\tau)$ is sated for the idempotent
class $\sfZ_0^{\bbR\bf{e}_1}\vee \sfZ_0^{\bbR\bf{e}_2}$ then the
factors
\[\xi_1 = \xi_2 := \zeta_0^{\tau^{\ \uhr\bbR\bf{e}_1}}\vee \zeta_0^{\tau^{\ \uhr\bbR\bf{e}_2}}\]
are characteristic for the above averages.
\end{prop}

\textbf{Proof of Theorem~\ref{thm:cts-poly-aves} from
Proposition~\ref{prop:char-cts-poly}}\quad Given
Proposition~\ref{prop:char-cts-poly} it suffices to consider the
averages $S_T(f_1,f_2)$ with each $f_i$ measurable with respect to
$\zeta_0^{\tau^{\ \uhr\bbR\bf{e}_1}}\vee \zeta_0^{\tau^{\
\uhr\bbR\bf{e}_2}}$.  By a simple approximation in $L^2(\mu)$ and
multilinearity, the convergence of these follows in turn if we know
it when $f_i = g_i\cdot h_i$ for some $g_1,g_2$ that are $\tau^{\
\uhr\bbR\bf{e}_1}$-invariant and $h_1,h_2$ that are $\tau^{\
\uhr\bbR\bf{e}_2}$-invariant.

Substituting this form into the definition of $S_T$, and using first
the invariance of $g_1$ and then that of $h_1$, we are left with the
averages
\begin{multline*}
S_T(f_1,f_2) = g_1\cdot \frac{1}{T}\int_0^T (h_1\circ
\tau^{t^2\bf{e}_1})((g_2\cdot h_2)\circ\tau^{t^2\bf{e}_1 +
t\bf{e}_2})\,\d t\\ = g_1\cdot \frac{1}{T}\int_0^T (h_1\cdot
g_2\cdot h_2)\circ\tau^{t^2\bf{e}_1 + t\bf{e}_2}\,\d t,
\end{multline*}
and these latter are now conventional polynomial ergodic averages,
which converge in $L^2(\mu)$ simply by spectral theory and the
corresponding result for the scalar averages $\frac{1}{T}\int_0^T
\exp(2\pi (at^2 + bt)\rm{i})\,\d t$, whose convergence follows from
the classical scalar-valued van der Corput estimate (or, indeed, in
the only nontrivial case $a\neq 0$, from a change of variables and
the classical evaluation of Fresnel integrals). \qed

\textbf{Remark}\quad Note that Proposition~\ref{prop:char-cts-poly}
is formulated in terms of satedness with respect to a single
idempotent class, rather than by appeal to a continuous analog of
the blanket notion of full isotropy satedness
(Definition~\ref{dfn:FIS}).  This is because the continuous group
$\bbR^2$ has uncountably many subgroups, and so we should need in
turn an analog of Theorem~\ref{thm:sateds-exist} that allows for
uncountably many idempotent classes.  This is impossible in general
without leaving the class of standard Borel spaces (although
presumably the class of `perfect' measure spaces, in the sense of
Section 451 of~\cite{FreVol4}, is still large enough), and so would
entail a barrage of new technical measure-theoretic details that we
prefer to avoid. \qed

\textbf{Proof of Proposition~\ref{prop:char-cts-poly}}\quad This
proof that starts with an important initial twist. We first note
that we may change variables in the integral
\[\frac{1}{T}\int_0^T(f_1\circ \tau^{t^2\bf{e}_1})(f_2\circ \tau^{t^2\bf{e}_1 + t\bf{e}_2})\,\d t\]
to $u := t^2$, and so obtain
\begin{eqnarray*}
S_{\sqrt{U}}(f_1,f_2) &=& \frac{1}{\sqrt{U}}\int_0^U (f_1\circ
\tau^{u\bf{e}_1})(f_2\circ \tau^{u\bf{e}_1 +
\sqrt{u}\bf{e}_2})\,\frac{\d u}{2\sqrt{u}}\\ &=&
\frac{1}{2}S'_U(f_1,f_2) + \frac{1}{\sqrt{U}}\int_0^U
\frac{1}{4}V^{-1/2}\cdot S'_V(f_1,f_2)\,\d V
\end{eqnarray*}
where
\[S'_U(f_1,f_2) := \frac{1}{U}\int_0^U(f_1\circ \tau^{u\bf{e}_1})(f_2\circ
\tau^{u\bf{e}_1 + \sqrt{u}\bf{e}_2})\,\d u.\] Thus, this change of
variables has revealed that the averages $S_T(f_1,f_2)$ are actually
`smoother' than the averages $S'_U(f_1,f_2)$, which involve only
linear and sub-linear exponents. In spite of of the non-integer
power $\sqrt{u}$ that has now appeared, we will now see that these
are quite simple for our purposes. (This crucial trick was pointed
out to me by Vitaly Bergelson.)

To complete the proof we show that \[S'_U(f_1,f_2)\not\to
0\quad\quad\Rightarrow\quad\quad \sfE_\mu(f_1\,|\,\zeta_0^{\tau^{\
\uhr\bbR\bf{e}_1}}\vee \zeta_0^{\tau^{\ \uhr\bbR\bf{e}_2}}) \neq
0.\] As usual, this begins with the van der Corput estimate (in its
version for continuous families of vectors, which is exactly
analogous to the discrete setting: see, for example, Section 1.9 of
Kuipers and Niederreiter~\cite{KuiNie74}), which after a little
re-arrangement gives that
\begin{eqnarray*}
&&S'_U(f_1,f_2)\not\to 0\\
\Rightarrow&& \frac{1}{H}\int_0^h\frac{1}{U}\int_0^U\int_X((f_1\circ
\tau^{h\bf{e}_1}\cdot \bar{f_1})\circ \tau^{u\bf{e}_1})\\
&&\quad\quad\quad\quad\quad\cdot ((f_2\circ \tau^{h\bf{e}_1 +
(\sqrt{u + h} - \sqrt{u})\bf{e}_2}\cdot \bar{f_2}) \circ
\tau^{u\bf{e}_1 + \sqrt{u}\bf{e}_2})\,\d\mu\,\d u\,\d h\\
&&= \frac{1}{H}\int_0^h\frac{1}{U}\int_0^U\int_X(f_1\circ
\tau^{h\bf{e}_1}\cdot \bar{f_1})\\
&&\quad\quad\quad\quad\quad\cdot ((f_2\circ \tau^{h\bf{e}_1 +
(\sqrt{u + h} - \sqrt{u})\bf{e}_2}\cdot \bar{f_2}) \circ
\tau^{\sqrt{u}\bf{e}_2})\,\d\mu\,\d u\,\d h \not\to 0\\
&&\quad\quad\quad\quad\quad\quad\quad\quad\quad\quad\quad\quad\hbox{as
$U\to\infty$ and then $H\to\infty$}.
\end{eqnarray*}

The important feature here is that for each fixed $h$ we have
\[\sqrt{u+h} - \sqrt{u}\to 0\quad\quad\hbox{as}\ u\to\infty,\]
and hence by the strong continuity of $\tau$ it follows that
\[\|f_2\circ \tau^{h\bf{e}_1 + (\sqrt{u + h} - \sqrt{u})\bf{e}_2} - f_2\circ \tau^{h\bf{e}_1}\|_2 \to 0\]
as $u \to\infty$.  From this it follows that for any fixed $h$ we
have
\begin{eqnarray*}
&&\frac{1}{U}\int_0^U\int_X(f_1\circ \tau^{h\bf{e}_1}\cdot
\bar{f_1})((f_2\circ \tau^{h\bf{e}_1 + (\sqrt{u + h} -
\sqrt{u})\bf{e}_2}\cdot \bar{f_2}) \circ
\tau^{\sqrt{u}\bf{e}_2})\,\d\mu\,\d u\\
&&\sim \frac{1}{U}\int_0^U\int_X(f_1\circ \tau^{h\bf{e}_1}\cdot
\bar{f_1})((f_2\circ \tau^{h\bf{e}_1}\cdot
\bar{f_2}) \circ \tau^{\sqrt{u}\bf{e}_2})\,\d\mu\,\d u\\
&&\to \int_X(f_1\circ \tau^{h\bf{e}_1}\cdot
\bar{f_1})\sfE_\mu(f_2\circ \tau^{h\bf{e}_1}\cdot
\bar{f_2}\,|\,\zeta_0^{\tau^{\ \uhr\bbR\bf{e}_2}})\,\d\mu\quad\quad\hbox{as}\ U\to\infty\\
&&= \int_{X^2}(f_1\circ \tau^{h\bf{e}_1}\cdot \bar{f_1})\otimes
(f_2\circ \tau^{h\bf{e}_1}\cdot
\bar{f_2})\,\d(\mu\otimes_{\zeta_0^{\tau^{\ \uhr\bbR\bf{e}_2}}}\mu)
\end{eqnarray*}
(this crucial simplification resulting from our change-of-variables
was pointed out to me by Vitaly Bergelson).  Now letting
$h\to\infty$ this simply tends to
\[\int_{X^2}(f_1\otimes f_2)\cdot g\,\d(\mu\otimes_{\zeta_0^{\tau^{\ \uhr\bbR\bf{e}_2}}}\mu)\]
for the $(\tau^{\otimes 2})^{\ \uhr\bbR\bf{e}_1}$-invariant function
\[g := \sfE_{\mu\otimes_{\zeta_0^{\tau^{\ \uhr\bbR\bf{e}_2}}}\mu}(f_1\otimes f_2\,|\,\zeta_0^{(\tau^{\otimes 2})^{\ \uhr\bbR\bf{e}_1}})\]

Hence, letting $(\t{X},\t{\mu}) = (X^2,\mu\otimes_{\zeta_0^{\tau^{\
\uhr\bbR\bf{e}_2}}}\mu)$, letting $\pi:\t{X} \to X$ and lifting
$\tau$ to the action $\t{\tau}$ defined by
\[\t{\tau}^{s\bf{e}_1 + t\bf{e}_2} := \tau^{s\bf{e}_1 + t\bf{e}_2}\times \tau^{s\bf{e}_1}\]
(noting that $\mu\otimes_{\zeta_0^{\tau^{\ \uhr\bbR\bf{e}_2}}}\mu$
is also invariant under the flow $t\mapsto \id_X\times
\tau^{t\bf{e}_2}$), we see that we have found an extension of
$(X,\mu,\tau)$ in which
\[\sfE_{\t{\mu}}(f_1\circ \pi\,|\,\zeta_0^{\tau^{\ \uhr\bbR\bf{e}_1}}\vee \zeta_0^{\tau^{\ \uhr\bbR\bf{e}_2}}) \neq 0,\]
and hence by satedness the analogous non-vanishing must have held
inside the original system $(X,\mu,\tau)$, as required. \qed

\textbf{Remark}\quad Although rather simple, it is worth noting that
the use of satedness still played a crucial r\^ole in the above
proof.  Without the initial assumption of satedness, the above
appeal to the van der Corput estimate combined with considerations
of the structure of the Furstenberg self-joining tell us that for
the averages $S_U(f_1,f_2)$ the pair of factors $\xi_1$, $\xi_2$ is
characteristic, where $\xi_i$ coordinatizes the maximal $\tau^{\
\uhr\bbR\bf{e}_i}$-isometric extension of the isotropy factor
$\zeta_0^{\tau^{\ \uhr\bbR(\bf{e}_1 - \bf{e}_2)}}$.  This much can
be argued following the same lines as Conze and Lesigne's initial
analysis in~\cite{ConLes84,ConLes88.1,ConLes88.2} of double linear
averages in discrete time for some system $T:\bbZ^2\curvearrowright
(X,\mu)$. Thus, allowing ourselves to assume that each $f_i$ is
$\xi_i$-measurable, it now follows that each $f_i$ may be
approximated by a function residing in a finite-rank $\tau^{\
\uhr\bbR\bf{e}_i}$-invariant module over the factor
$\zeta_0^{\tau^{\ \uhr\bbR(\bf{e}_1 - \bf{e}_2)}}$.  In Conze and
Lesigne's setting (with $T$ and $\bbZ$ in place of $\tau$ and
$\bbR$) this leads directly to a proof of convergence, because when
written in terms of unitary cocycles describing these finite-rank
modules the double linear averages become simply averages for some
new `combined' finite rank module over the single system
$T|_{\zeta_0^{T^{\ \uhr\bbZ(\bf{e}_1 - \bf{e}_2)}}}$, to which the
usual mean ergodic theorem can be applied.  However, in our setting
matters are not so simple, since even after approximating and then
using a representation in terms of unitary cocycles in this way, the
expression that results still involves two different polynomials in
the exponents, and so it is not clear how to realize it as some kind
of more classic ergodic average. Once we assumed satedness, this
problem vanished because the structure of a finite-rank $\tau^{\
\uhr\bbR\bf{e}_i}$-invariant module over $\zeta_0^{\tau^{\
\uhr\bbR(\bf{e}_1 - \bf{e}_2)}}$ is replaced by that of the factor
$\zeta_0^{\tau^{\ \uhr\bbR\bf{e}_i}}\vee \zeta_0^{\tau^{\
\uhr\bbR(\bf{e}_1 - \bf{e}_2)}}$, for which more explicit
simplifications to our averages are possible, as exhibited above.
\fin

\textbf{Remark}\quad Since the present paper was first submitted, the above argument has been generalized in~\cite{Aus--normconvctstime} to handle all continuous-time polynomial nonconventional ergodic averages for actions of $\bbR^d$.  A slightly different approach, still using satedness and giving a further generalization to actions of nilpotent Lie groups, was then presented in~\cite{Aus--joiningequidistnilpotent}. \fin

\section{Next steps}

This paper has begun to showcase the far-reaching consequences of satedness in the study of nonconventional
ergodic averages, but its larger purpose is to prepare the ground
for its sequel~\cite{Aus--lindeppleasant2}.  There we will turn to
nonconventional averages which require rather more elaborate new
arguments. Indeed, only after several more technical steps will we
be able to address even one new case of convergence for polynomial
averages in discrete time: that of
\[\frac{1}{N}\sum_{n=1}^N (f_1\circ T_1^{n^2})(f_2\circ T_1^{n^2}T_2^n)\quad\quad\hbox{for}\ (T_1,T_2):\bbZ^2\curvearrowright (X,\mu),\ f_1,f_2 \in L^\infty(\mu)\]
(as in Theorem~\ref{thm:polyconv} above).

The analysis of these will rely heavily on some auxiliary results concerning the triple linear averages
\[\frac{1}{N}\sum_{n=1}^N (f_1\circ T^{n\bf{p}_1})(f_2\circ T^{n\bf{p}_2})(f_3\circ T^{n\bf{p}_3})\]
for some action $T:\bbZ^2\curvearrowright (X,\mu)$ and three
directions $\bf{p}_1$, $\bf{p}_2$, $\bf{p}_3 \in \bbZ^2$ enjoying
some linear dependence.  Such linear averages arise naturally from
the above quadratic averages upon a single application of the van
der Corput estimate, and so we will construct a useful notion of
pleasant extension for the quadratic averages by first developing
such a notion for these triple linear averages and then showing how
the resulting characteristic factors can be simplified further.

The point is that, although convergence is known for triple linear
averages such as the above, the use of extensions to prove this
in~\cite{Aus--nonconv} forgets the
linear dependence of the $\bf{p}_i$, effectively replacing the
$T^{\bf{p}_i}$ with three independent commuting transformations on
the extended system. We cannot afford this freedom in the study of
the quadratic averages, because after passing to such a
$\bbZ^3$-system it is not clear how the quadratic averages of
interest can even be sensibly interpreted.  Our main task,
therefore, will be to see how simple a triple of characteristic
factors can be obtained for the above linear nonconventional
averages while preserving the algebraic relations of the original
$\bbZ^2$-action. It will turn out that we can do quite well, except that in addition to the factors that contribute to each $\xi_i$ in
Theorem~\ref{thm:char-retaining-roots} we must now involve some systems on which our $\bbZ^2$-action is by commuting rotations on a two-step nilmanifold, thus re-establishing contact with earlier works such as~\cite{HosKra05,Zie07} and their forerunners in which the relevance of these was made clear in the setting of $\bbZ$-actions. 

\parskip 0pt

\bibliographystyle{abbrv}
\bibliography{bibfile}

\vspace{10pt}

\small{\textsc{Courant Institute, New York University, New York, NY 10012, USA}}

\vspace{5pt}

\small{Email: \verb|tim@cims.nyu.edu|}

\vspace{5pt}

\small{URL: \verb|http://www.cims.nyu.edu/~tim|}

\end{document}